\documentclass[12pt,a4paper]{amsart} 
\usepackage[utf8]{inputenc}
\usepackage[english]{babel}

\usepackage{amsmath}
\usepackage{amsfonts}
\usepackage{amssymb}
\usepackage{color}
\usepackage{ifthen}

\usepackage{amsthm}

\usepackage{graphics}
\usepackage{epsfig}

\usepackage{url}
\usepackage{hyperref}

\usepackage[a4paper,top=2.5cm,bottom=2.8cm,
left=3.2cm,right=3.2cm]{geometry}
\markleft{\hfill \textsc{Computer Investigation of Properties of the Gergonne Point of a Triangle} \hfill}
\newcommand{\degrees}{^\circ}
\newcommand{\GE}{GeometricExplorer}

\newcommand{\citeproblem}[1]{\textcolor{blue}{\cite{#1}}}

\newcounter{propertyNumber}

\newsavebox{\figbox}
\newlength{\figheight}

\long\def\void#1{}

\newcommand{\prop}[7]{
\medskip
\smallskip
\goodbreak
\begin{raggedright}
\addtocounter{propertyNumber}{1}%
\savebox{\figbox}{\includegraphics[width=0.4\linewidth]{#6.png}}%
\textbf{Property \thesubsection.\thepropertyNumber.} \textcolor{blue}{\textbf{#1}}
\settoheight{\figheight}{\includegraphics[width=#7\linewidth]{#6.png}}%
\fbox{\parbox[b][\figheight]{0.57\linewidth}{
\raggedright
\textbf{Start with:} #2\\
\ifthenelse{\equal{#3}{}}{}
{\medskip
\textbf{Step 1:} #3\\
}
\medskip
\textbf{Conclusion:} #4
\vfill
\textbf{References:} #5}
\vrule width 1pt
\usebox{\figbox}}
\\
\end{raggedright}
}

\newcommand{\resetPropertyNumber}{\setcounter{propertyNumber}{0}}
\begin{document}
% ===================
International Journal of  Computer Discovered Mathematics (IJCDM) \\
ISSN 2367-7775 \copyright IJCDM \\
Volume 6, 2021, pp.X--X  \\
Received XX January 2021. Published on-line XX XXX 2021 \\ 
web: \url{http://www.journal-1.eu/} \\
\copyright The Author(s) This article is published 
with open access\footnote{This article is distributed under the terms of the Creative Commons Attribution License which permits any use, distribution, and reproduction in any medium, provided the original author(s) and the source are credited.}. \\
% ===========================   
\bigskip
\medskip

\begin{center}
	{\Large \textbf{Computer Investigation of Properties}} \\
	\medskip
	{\Large \textbf{of the Gergonne Point of a Triangle}} \\
	\bigskip

	\textsc{Stanley Rabinowitz$^a$ and Ercole Suppa$^b$} \\

	$^a$ 545 Elm St Unit 1,  Milford, New Hampshire 03055, USA \\
	e-mail: \href{mailto:stan.rabinowitz@comcast.net}{stan.rabinowitz@comcast.net}\footnote{Corresponding author} \\
	web: \url{http://www.StanleyRabinowitz.com/} \\
	
	$^b$ Via B. Croce 54, 64100 Teramo, Italia \\
	e-mail: \href{mailto:ercolesuppa@gmail.com}{ercolesuppa@gmail.com} \\
     %draft: Dec. 30, 2020
\end{center}
\bigskip

% ==============================
\textbf{Abstract.} The incircle of a triangle touches the sides of the triangle in three points.
It is well-known that the lines from these points to the opposite vertices meet
at a point known as the Gergonne point of the triangle.
We use a computer to discover and catalog properties of the Gergonne point.



\newcommand{\Gergonne}{\mathtt{Gergonne}}
\newcommand{\Nagel}{\mathtt{Nagel}}
\newcommand{\Spieker}{\mathtt{Spieker}}
\newcommand{\orthocenter}{\mathtt{orthocenter}}
\newcommand{\incenter}{\mathtt{incenter}}
\newcommand{\centroid}{\mathtt{centroid}}
\newcommand{\para}{\mathtt{parallel}}
\newcommand{\perpen}{\mathtt{perpendicular}}
\newcommand{\incircle}{\mathtt{incircle}}
\newcommand{\midpoint}{\mathtt{midpoint}}
\newcommand{\mittenpunkt}{\mathtt{mittenpunkt}}
\newcommand{\Feuerbach}{\mathtt{Feuerbach}}
\newcommand{\ninepointcenter}{\mathtt{ninepointcenter}}
\newcommand{\colline}{\mathtt{colline}}
\newcommand{\concur}{\mathtt{concur}}
\newcommand{\foot}{\mathtt{foot}}
\newcommand{\noproof}{\textcolor{red}{*}}

% ================================
% 1 Introduction 
% ================================
%\tableofcontents
\section{Introduction}
\label{section:intro}

The line from a vertex of a triangle to the contact point of the incircle with
the opposite side is called a \textit{Gergonne cevian}.
It is well known that the three Gergonne cevians of a triangle
meet in a point. This point is known as the Gergonne point of the triangle.
It is shown in the following illustration.
%It is shown in Figure \ref{fig:GergonnePoint}.

\begin{figure}[h!t]
\centering
\includegraphics[width=0.4\linewidth]{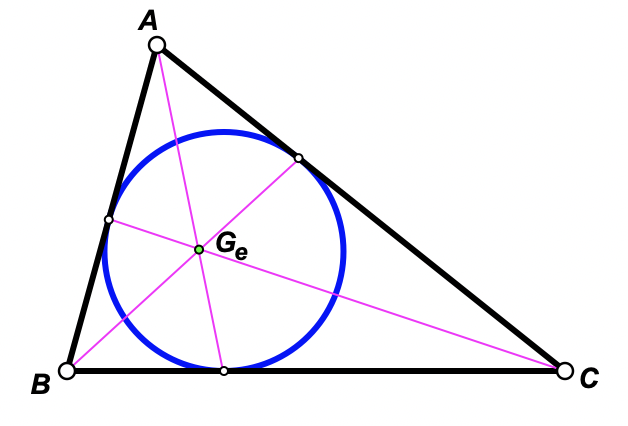}
%\caption{Gergonne point}
\label{fig:GergonnePoint}
\end{figure}

We use a computer to investigate properties of the Gergonne point.
Some of the properties found have been discovered before, and we will give a reference,
if known, to these results. In other cases, the property is new.

We systematically search for properties of the Gergonne point by starting with a
configuration that consists of a triangle and its Gergonne point.
Then we apply a sequence of 0 to 2 basic
constructions to this configuration.
Basic constructions consist of typical ruler and compass constructions such as finding
the intersection of a line with a line or a circle, dropping a perpendicular from a
point to a line, drawing parallel and perpendicular lines, constructing an incircle or circumcircle, etc. We also included constructions involving triangle centers, such as constructing centroids,
orthocenters, and symmedian points of a triangle.
Some Apollonius circle constructions, including center and all touch points
(PPP, PPL, PPC, LLL, LLP, LLC, CLP) were also included.
This includes drawing incircles, excircles, circumcircles, and mixtilinear incircles.
All possible sequences of such constructions are generated.

\bigskip
\begin{center}
\begin{tabular}{|l|l|}
\hline
\multicolumn{2}{|c|}{\large \strut Notation for Constructions}\\
\hline
\textbf{Notation}&\textbf{Description}\\
\hline
$\triangle ABC$&Triangle $ABC$\\
$a$, $b$, $c$&The lengths of the sides of $\triangle ABC$\\
$s$&$(a+b+c)/2$\\
$\odot ABC$&The circle through points $A$, $B$, and $C$\\
$r(ABC)$&The radius of the incircle of $\triangle ABC$\\
$R(ABC)$&The radius of the circumcircle of $\triangle ABC$\\
$[F]$&The area of figure $F$\\
$\Gergonne[ABC]$&The Gergonne point of $\triangle ABC$\\
$\Nagel[ABC]$&The Nagel point of $\triangle ABC$\\
$\Spieker[ABC]$&The Spieker center of $\triangle ABC$\\
$\orthocenter[ABC]$&The orthocenter of $\triangle ABC$\\
$\mittenpunkt[ABC]$&The mittenpunkt of $\triangle ABC$\\
$\Feuerbach[ABC]$&The Feuerbach point of $\triangle ABC$\\
$\para[P,AB]$&The line through $P$ parallel to $AB$\\
$\perpen[P,AB]$&The line through $P$ perpendicular to $AB$\\
$\foot[A,BC]$&The foot of the perpendicular from $A$ to $BC$\\
$\incircle[ABC]$&The incircle of $\triangle ABC$\\
$\incenter[ABC]$&The incenter of $\triangle ABC$\\
$\centroid[ABC]$&The centroid of $\triangle ABC$\\
$\midpoint[AB]$&The midpoint of line segment $AB$\\
$\ninepointcenter[ABC]$&The nine-point center of $\triangle ABC$\\
$\colline[A,B,C]$&The points $A$, $B$, and $C$ colline.\\
$\concur[AB,CD,EF]$&The lines $AB$, $CD$, and $EF$ are concurrent.\\
\hline
\end{tabular}
\end{center}
\bigskip

After each construction is applied, the resulting figure is analyzed for basic properties, such
as three points being collinear, three lines being concurrent, checking for perpendicular or parallel lines,
congruent incircles, linear and quadratic relationships between lengths of line segments, etc.
Trivial properties are excluded, as well as standard properties of
the construction that is applied. The process is carried out by a program called \GE\
that we wrote specifically for this purpose.
In some cases, when GeometricExplorer found interesting numerical properties, we used
a Mathematica program to find a more general result.

\textbf{Convention:} Gergonne points in our figures will be colored green.

\textbf{Proofs:}
The program we used, \GE, works by starting with numerical figures and examining the figures constructed
using 15 decimal places of accuracy. Although any result found is most likely true, this
does not constitute a proof. Using barycentric coordinates and Mathematica, we have found proofs for all of the properties discovered. The proofs are in the Mathematica notebook supplied with
the supplementary material associated with this paper.

\textbf{References:}
A reference in blue indicates that the result was discovered by our investigation and is believed to be new and was posted to social media in the hopes that someone would find a purely geometric solution. 
An asterisk after a reference indicates that no one has found a purely geometrical proof as of January 2021.

\raggedbottom

%******************
%
%     Triangle and Gergonne Point
%
%******************

%\newpage
\section{A Triangle and its Gergonne Point}

\resetPropertyNumber

In this section, we report on properties found for a starting figure consisting of a
triangle together with its Gergonne point.

\subsection{Arbitrary Triangle}\ \medskip

\textbf{Intrinsic Properties}\ \medskip

An \textit{intrinsic property} of a figure is a property of that figure
before any constructions are applied. For example, an intrinsic property of an
equilateral triangle is that the angles all have measure $60\degrees$.

\prop{(Spoke Formula)}
{Triangle $ABC$.\\
$D=\Gergonne[ABC]$. $x=DA$.}
{}
{\\
\vspace*{4pt}
$\displaystyle x=\sqrt\frac{a(b+c-a)^3\left(a(a+b+c)-2(b-c)^2\right)}{(a^2+b^2+c^2-2ab-2bc-2ca)^2}$.}
{See supplementary notebook.}
{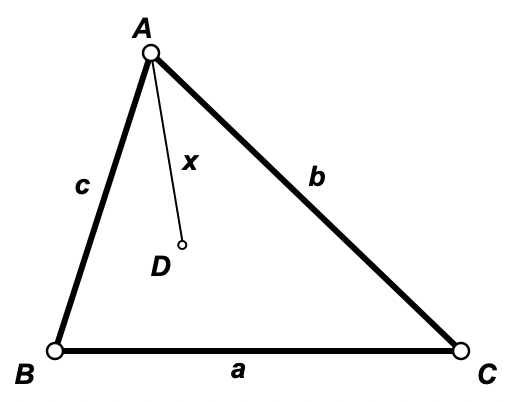}{0.4}

\prop{(Tripolar Coordinates)}
{Triangle $ABC$.\\
$D=\Gergonne[ABC]$. $x=DA$, $y=DB$.}
{}
{\\
\vspace*{8pt}
$\frac{x}{y}=\frac{b+c-a}{c+a-b}\sqrt{\frac{a(b+c-a)(a^2+ab-2b^2+ac+4bc-2c^2)}{b(c+a-b)(b^2+bc-2c^2+ba+4ca-2a^2)}}$.
}
{This follows from the SpokeFormula. See also \cite{ETC7}.}
{tripolar}{0.4}

\prop{(Barycentric Coordinates)}
{Triangle $ABC$.\\
$D=\Gergonne[ABC]$.}
{}
{\\
\vspace*{14pt}
$\displaystyle\frac{[BDC]}{[CDA]}=\frac{c+a-b}{b+c-a}$.
}
{\cite{ETC7}}
{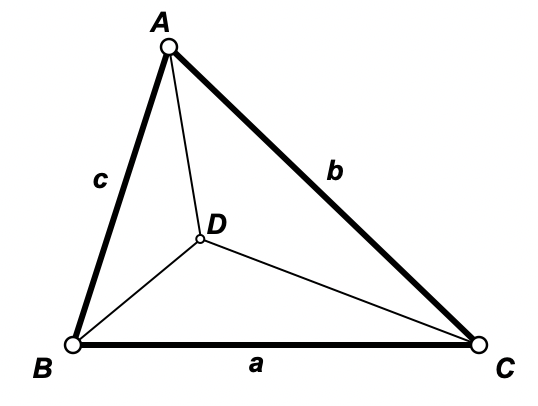}{0.4}

\prop{(Area Formula)}
{Triangle $ABC$.\\
$D=\Gergonne[ABC]$.}
{}
{\\
\vspace*{14pt}
$\displaystyle [BCD]=\frac{(a+b-c)(a-b+c)}{2ab+2bc+2ca-a^2-b^2-c^2}K$,\\
\vspace*{14pt}
where $K=\sqrt{s(s-a)(s-b)(s-c)}$.
}
{See supplementary notebook.}
{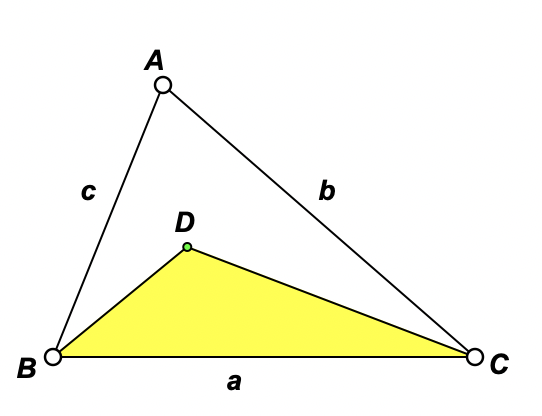}{0.43}

\bigskip
\textbf{1-step Properties}\ \medskip

A \textit{1-step property} of a figure is a property found
after applying exactly one basic construction to the given figure.

\prop{(Defining Property of the Gergonne Point)}
{Triangle $ABC$.\\$D=\Gergonne[ABC]$.}
{$E=\incircle[ABC]\cap BC$}
{$\colline[A,D,E]$.}
{\cite[p.~160]{Altshiller-Court}, \cite[p.~184]{Johnson}}
{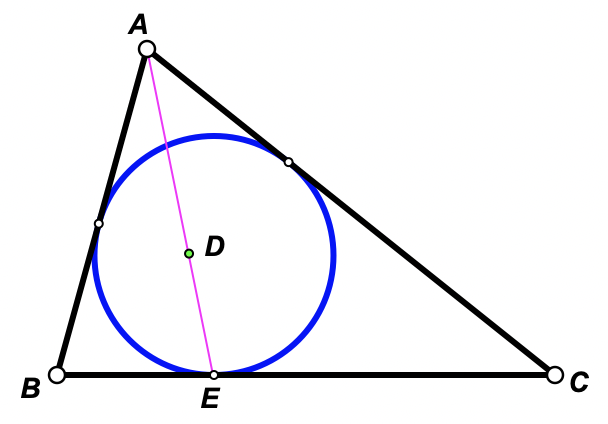}{0.4}

\prop{(LLP Property)}
{Triangle $ABC$.\\
$D=\Gergonne[ABC]$.}
{A circle, center $E$, passes through $D$ and is tangent to $AB$ and $AC$.}
{$DE\perp BC$.}
{\citeproblem{RG5901}}
{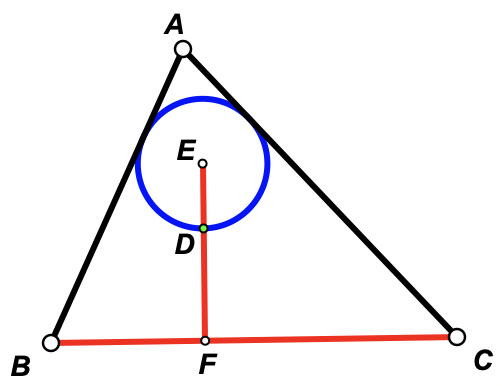}{0.4}

\prop{(Intouch Symmedian Point)}
{Triangle $ABC$.\\
$G_e=\Gergonne[ABC]$.}
{The incircle of $\triangle ABC$ touches $BC$ at $D$, touches $AC$ at $E$, and touches $AB$ at $F$.}
{$G_e$ is the symmedian point of $\triangle DEF$.}
{\cite{ETC7}}
{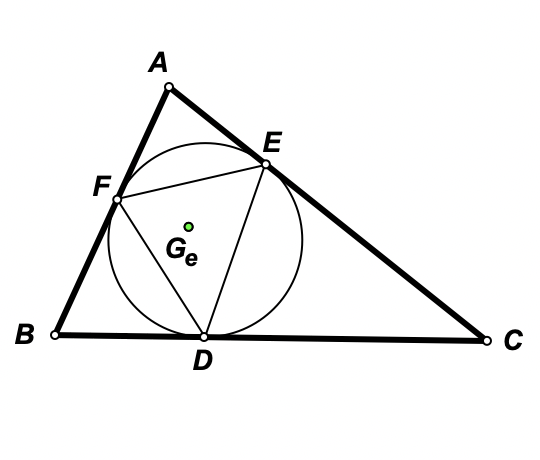}{0.4}

\prop{}
{Triangle $ABC$.\\
$D=\Gergonne[ABC]$.}
{The incircle of $\triangle ABC$ touches $BC$ at $E$, touches $AC$ at $F$, and touches $AB$ at $G$.}
{\\
\vspace*{5pt}
$R(AFD)+R(BED)=R(ABD)$.}
{\citeproblem{RG6485}}
{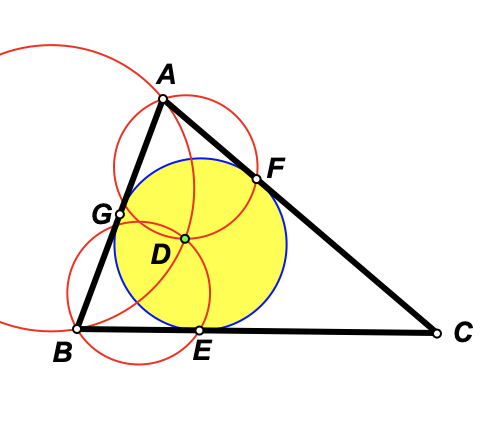}{0.4}

\bigskip
\subsection{Triangles with Angle Restrictions}\ \medskip

\resetPropertyNumber

In this section, we report on properties found when the starting triangle
has restrictions placed on its angles.

\medskip
\textbf{Intrinsic Properties}\ \medskip

No intrinsic properties were found.

\bigskip
\textbf{1-step Properties}\ \medskip

\prop{}
{Triangle $ABC$.
\\$D=\Gergonne[ABC]$.
}
{A circle, center $E$, passes through $D$ and is tangent to $AB$ at $F$
and tangent to $BC$.}
{$\angle DEF+\angle BAC=180\degrees$.}
{\cite{vanLamoen5985}}
{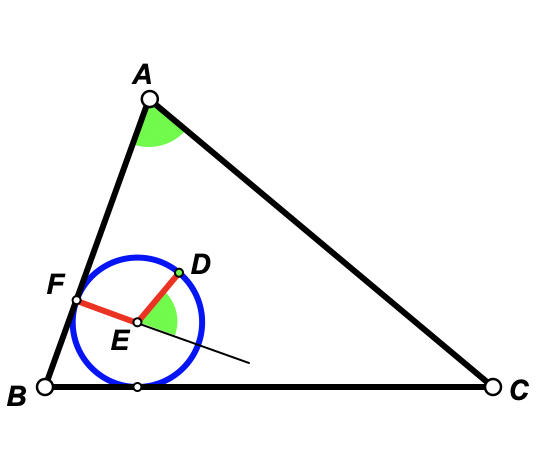}{0.4}

\prop{}
{Triangle $ABC$.\\
$\triangle ABC$ satisfies $\angle A=120\degrees$.\\
$D=\Gergonne[ABC]$.
}
{A circle, center $E$, passes through $D$ and is tangent to $AB$ at $F$
and tangent to $BC$.}
{$\triangle DEF$ is equilateral.}
{\citeproblem{RG5985}, special case of previous property}
{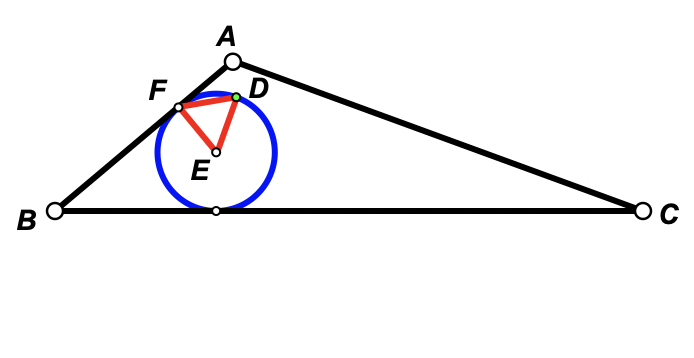}{0.65}

\prop{}
{Right Triangle $ABC$.\\$D=\Gergonne[ABC]$.\\
$\triangle ABC$ satisfies $\angle A=75\degrees$ and $\angle C=15\degrees$.
}
{$H=\orthocenter[ACD]$.}
{$[\odot ABC]+[\odot ACD]=[\odot BDH]$.}
{\citeproblem{RG5984}\noproof}
{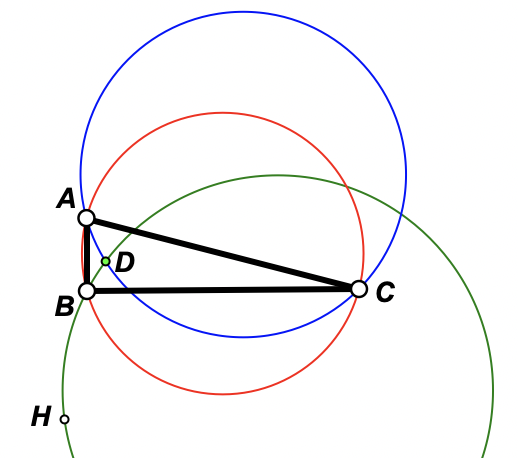}{0.4}

\prop{}
{Right triangle $ACB$.\\$D=\Gergonne[ABC]$.\\
$\triangle ABC$ satisfies $\angle C=90\degrees$.
}
{$E$ is the center of a circle that touches $BC$ at $C$ and touches $AB$ at $G$.}
{$[BDE]=[CDG]$.}
{\citeproblem{RG6002}\noproof}
{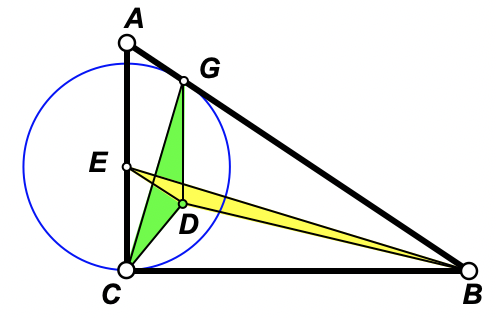}{0.4}

\prop{}
{Right triangle $ABC$.\\
$\triangle ABC$ satisfies $\angle B=90\degrees$.\\
$D=\Gergonne[ABC]$.
}
{The $B$-mixtilinear incircle touches $\odot ABC$ at $E$.}
{$\angle DBC=\angle DEC$.}
{\citeproblem{RG6026}\noproof}
{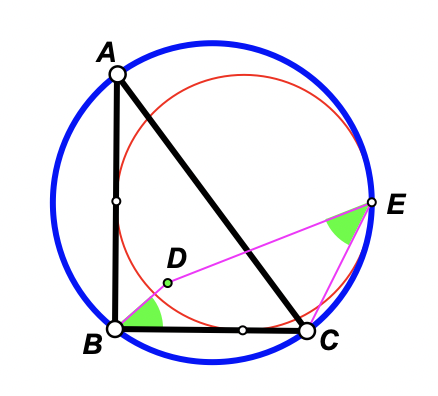}{0.4}

\newpage
\subsection{Triangles with Side Restrictions}\ \medskip

\resetPropertyNumber

In this section, we report on properties found when the starting triangle
has restrictions placed on its sides.

\bigskip
\textbf{Intrinsic Properties}\ \medskip

\prop{}
{Triangle $ABC$.\\
$\triangle ABC$ satisfies $a:b:c=7:9:10$.\\
$D=\Gergonne[ABC]$.
}
{}
{$AD=2CD$.}
{This follows from the Tripolar Coordinates.}
{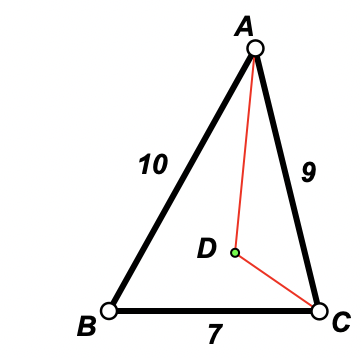}{0.4}

\prop{}
{Triangle $ABC$.\\
$\triangle ABC$ satisfies $a:b:c=5:8:9$.\\
$D=\Gergonne[ABC]$.
}
{}
{$AD=3CD$.}
{This follows from the Tripolar Coordinates.}
{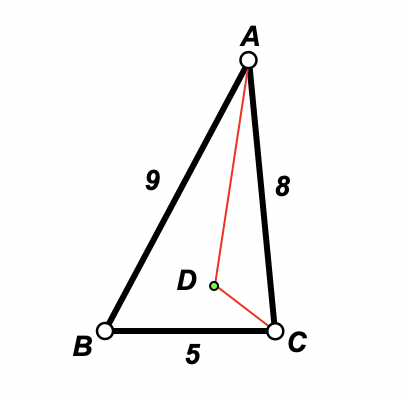}{0.4}

\prop{}
{Triangle $ABC$.\\
$\triangle ABC$ satisfies $a:b:c=4:9:10$.\\
$D=\Gergonne[ABC]$.
}
{}
{$AD=5CD$.}
{This follows from the Tripolar Coordinates.}
{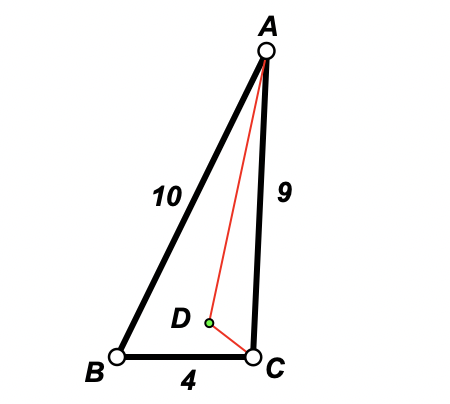}{0.4}

\prop{}
{Triangle $ABC$.\\
$\triangle ABC$ satisfies $as=b^2+bc+c^2$.\\
$D=\Gergonne[ABC]$.
}
{}
{The $A$-mixtilinear incircle of $\triangle ABC$ passes through $D$.
}
{\cite{Capitan6450}}
{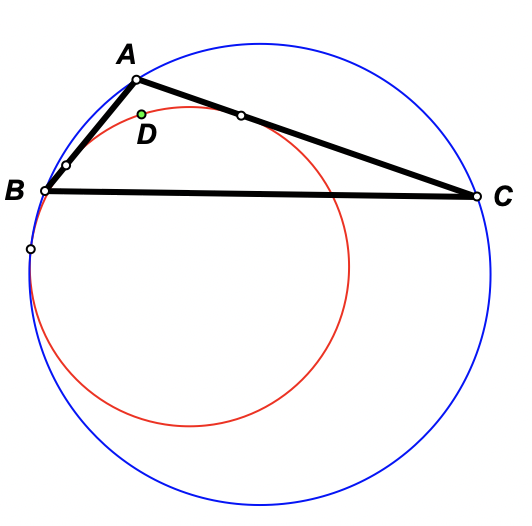}{0.43}

\bigskip
\bigskip
\textbf{1-step Properties}\ \medskip

\prop{}
{Triangle $ABC$.\\
$\triangle ABC$ satisfies $b=2a$.\\
$D=\Gergonne[ABC]$.
}
{The $C$-mixtilinear incircle touches $AC$ at $F$ and touches $\odot ABC$ at $E$.}
{$BD\parallel EF$.}
{\citeproblem{RG6005}\noproof}
{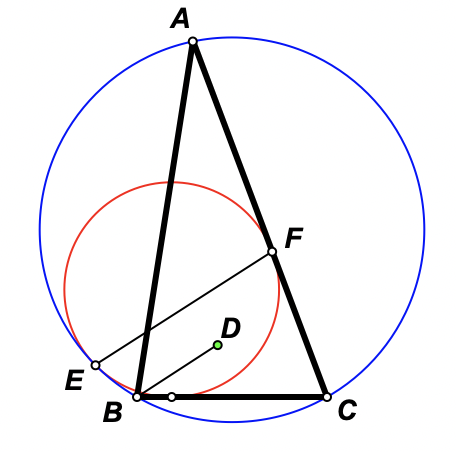}{0.4}

\prop{}
{Triangle $ABC$.\\
$\triangle ABC$ satisfies $a=3(b-c)$.\\
$D=\Gergonne[ABC]$.
}
{$E=\midpoint[AC]$.}
{$[ABD]=[CED]$.}
{\citeproblem{RG6008}}
{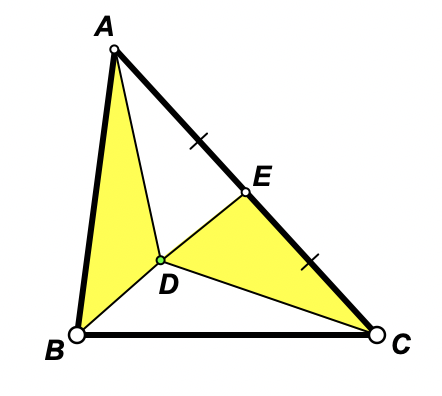}{0.4}

\prop{}
{Right triangle $ABC$.\\
$\triangle ABC$ satisfies $\angle C=90\degrees$.\\
$D=\Gergonne[ABC]$.
}
{$E=\perpen[D,BD]\cap BC$}
{$\displaystyle\frac{BD}{DE}=\sqrt{\frac{2(b+c)}{c-a}}$.}
{\cite{Suppa6012}}
{345perpa}{0.4}

\prop{}
{Right triangle $ABC$.\\
$\triangle ABC$ satisfies $\angle C=90\degrees$.\\
$D=\Gergonne[ABC]$.
}
{$E=\perpen[D,AD]\cap AC$}
{$\displaystyle\frac{AD}{DE}=\sqrt{\frac{2bs}{(c-b)(s-a)}}$.}
{\cite{Suppa6012}}
{345perpb}{0.4}

\prop{}
{Isosceles triangle $ABC$.\\
$\triangle ABC$ satisfies $a=c$.\\
$D=\Gergonne[ABC]$.}
{$E=\para[A,BC]\cap CD$.}
{\\
\vspace*{8pt}
$\displaystyle\frac{1}{BC}+\frac{1}{AE}=\frac{2}{AC}$.}
{\citeproblem{Suppa4247}}
{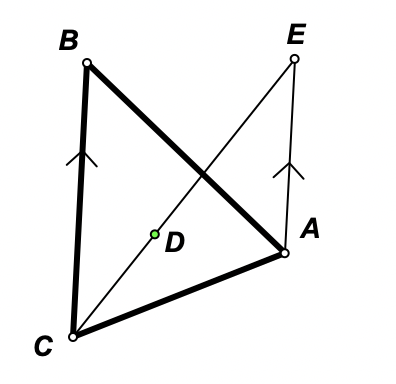}{0.4}

%\newpage
\bigskip
\subsection{A Triangle with a Gergonne Point and Another Center}\ \medskip

\resetPropertyNumber

In this section, we report on properties found when including the Gergonne point and another triangle center in the starting triangle.

\newpage
\prop{}
{Triangle $ABC$.\\
$\triangle ABC$ satisfies $\angle B=60\degrees$.\\
$D=\Gergonne[ABC]$.\\
$E=\incenter[ABC]$.
}
{}
{$\angle ABD=\angle AED$.}
{\citeproblem{RG7070}}
{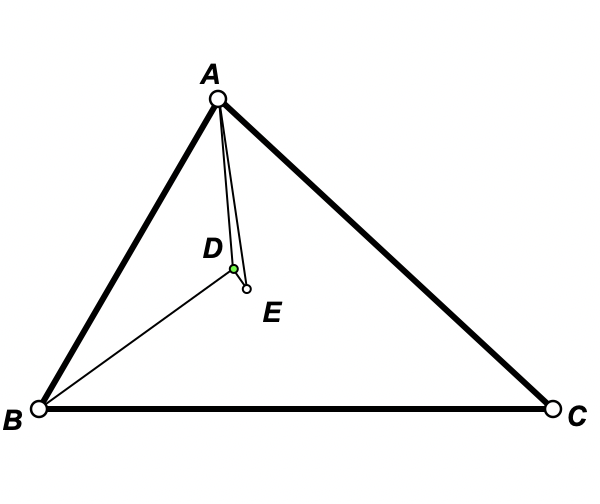}{0.4}

\prop{}
{Triangle $ABC$.\\
$D=\Gergonne[ABC]$.\\
$E=\mittenpunkt[ABC]$.
}
{}
{$[ABDE]=[BCE]$.}
{\citeproblem{RG7067}}
{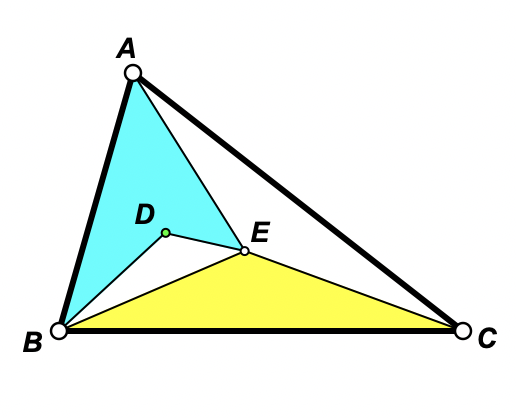}{0.4}

\prop{}
{Triangle $ABC$.\\
$\triangle ABC$ satisfies $2a=b+c$.\\
$D=\Gergonne[ABC]$.\\
$E=\Nagel[ABC]$.
}
{}
{$\angle BAD=\angle DEG$.}
{\citeproblem{RG7065}}
{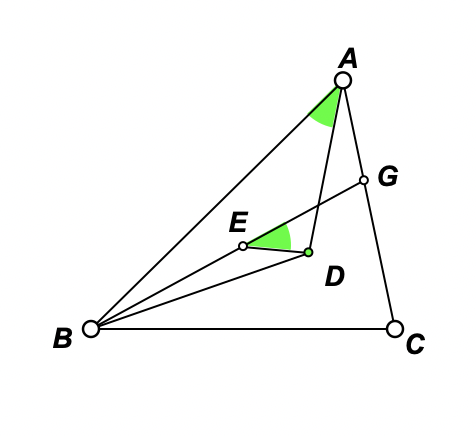}{0.4}

\prop{}
{Triangle $ABC$.\\
$\triangle ABC$ satisfies $2(s-a)(s-c)=b(s-b)$.\\
$D=\Gergonne[ABC]$.\\
$E=\Feuerbach[ABC]$.
}
{}
{$B$, $E$, and $D$ colline and
$BE=DE$.}
{\citeproblem{Suppa7073}}
{GEX11-9-8-5colline}{0.4}

\prop{}
{Triangle $ABC$.\\
$\triangle ABC$ satisfies $(s-b)(s-c)=a(s-a)$.\\
$D=\Gergonne[ABC]$.\\
$E=\Feuerbach[ABC]$.
}
{}
{$AE=2DE$.}
{\citeproblem{Suppa7073}}
{GEX11-9-8-5double}{0.4}

\prop{}
{Triangle $ABC$.\\
$\triangle ABC$ satisfies $a:b:c=9:7:4$.\\
$D=\Gergonne[ABC]$.\\
$E=\Nagel[ABC]$.
}
{}
{$CE=DE$.}
{\citeproblem{RG7074}}
{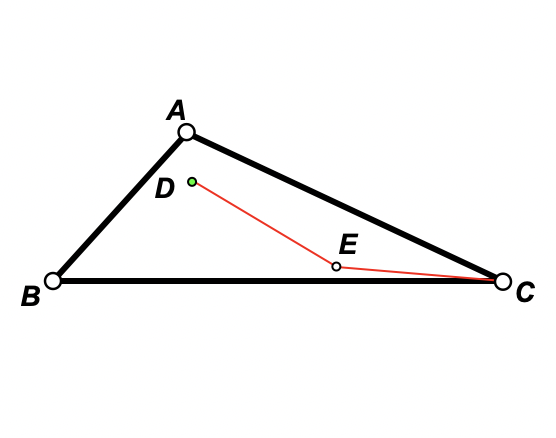}{0.4}
%2 a^3-2 a^2 b-2 a b^2+2 b^3-5 a^2 c+6 a b c-5 b^2 c+4 a c^2+4 b c^2-c^3 = 0

\prop{}
{Triangle $ABC$.\\
$\triangle ABC$ satisfies $a:b:c=9:7:4$.\\
$D=\Gergonne[ABC]$.\\
$E=\centroid[ABC]$.
$F=\mittenpunkt[ABC]$.
}
{}
{$AB\perp DE$ and $AD\perp BF$.}
{\citeproblem{RG7076}}
{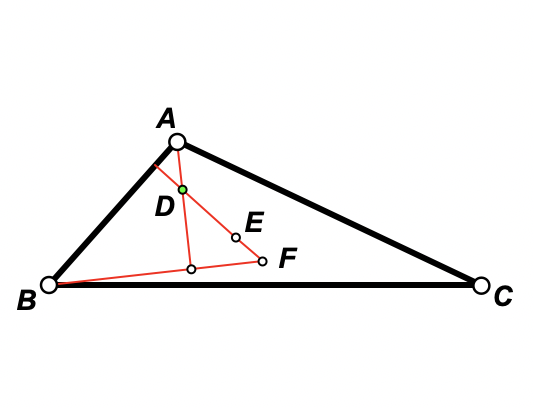}{0.4}

\bigskip
\textbf{Note:} Our investigation found many other perpendicularities involving
the Gergonne point and another triangle center in various shape triangles.
We omit them from this report since there were so many instances.

%******************
%
%     Triangle with Gergonne cevian
%
%******************

\newpage
\section{A Triangle with a Gergonne Cevian}

\resetPropertyNumber

In this section, we report on properties found for a starting figure consisting of a
triangle and one of its Gergonne cevians.

\subsection{Intrinsic Properties}\ \medskip

\prop{(Trace Lengths)}
{\\Triangle $ABC$ with Gergonne cevian $AE$.\\
$D=\Gergonne[ABC]$.
}
{}
{\\
\vspace*{6pt}
$BE=s-b$;\qquad
$CE=s-c$;\\
\vspace*{6pt}
$\displaystyle\frac{BE}{CE}=\frac{s-b}{s-c}$.}
{\cite[p.~87]{Altshiller-Court}}
{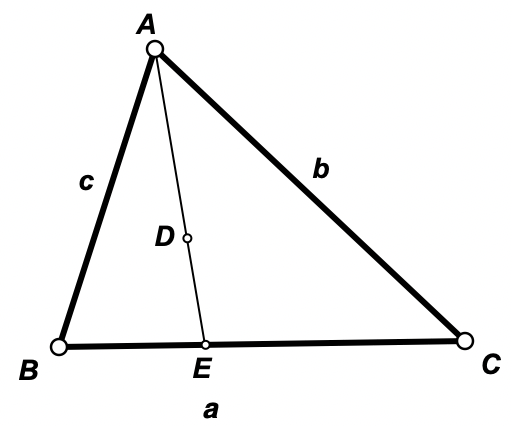}{0.4}

\prop{(Cevian Division)}
{\\Triangle $ABC$ with Gergonne cevian $AE$.\\
$D=\Gergonne[ABC]$.
}
{}
{
$\displaystyle\frac{AD}{DE}=\frac{a(s-a)}{(s-b)(s-c)}$.}
{\cite[p.~164]{Altshiller-Court}}
{traceRatio}{0.4}

\prop{(Cevian Length)}
{\\Triangle $ABC$ with Gergonne cevian $AE$.\\
$D=\Gergonne[ABC]$.
}
{}
{\\
\vspace*{4pt}
$\displaystyle AE^2=\frac{(s-a)(as-(b-c)^2)}{a}$.}
{See supplementary notebook.}
{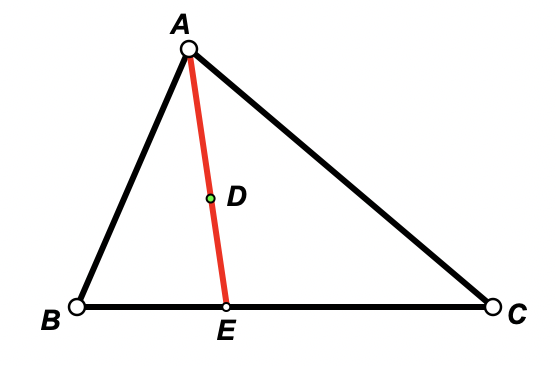}{0.4}
% was \frac{(b+c-a)(a^2-2b^2-2c^2+ab+4bc+ac)}{4a}

\newpage
\prop{(Split Perimeter Property)}
{\\Triangle $ABC$ with Gergonne cevian $AE$.\\
$D=\Gergonne[ABC]$.
}
{}
{$AB+CE=AC+BE=s$.}
{This follows from the trace lengths.}
{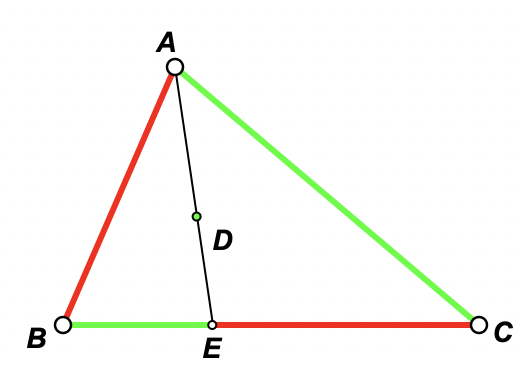}{0.4}

\prop{}
{\\Triangle $ABC$ with Gergonne cevian $BD$.\\
$\triangle ABC$ satisfies $\angle B=90\degrees$.
}
{}
{$AD\cdot CD=[ABC]$.}
{\citeproblem{RG6084}}
{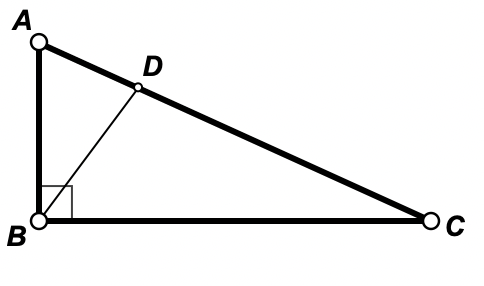}{0.4}

\prop{(Kissing Circles Theorem)}
{\\Triangle $ABC$ with Gergonne cevian $AD$.
}
{}
{\\Incircles of $\triangle ABD$ and $\triangle ADC$ touch.}
{\cite{kov}
%This is a special case of a theorem about tangential quadrilaterals.
%Note that $ABDC$ can be considered a tangential quadrilateral.
}
{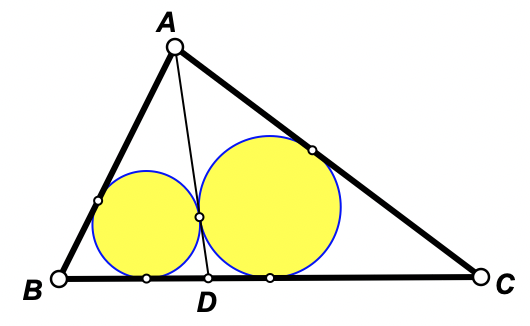}{0.4}

\prop{}
{\\Triangle $ABC$ with Gergonne cevian $CD$.\\
$\triangle ABC$ satisfies $3a=b+c$.
}
{}
{$AD=BC$.}
{This follows from the trace lengths.}
{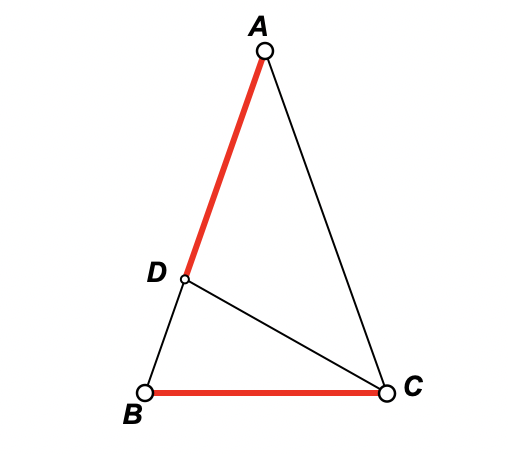}{0.4}

\newpage
\prop{}
{\\Triangle $ABC$ with Gergonne cevian $CD$.\\
$\triangle ABC$ satisfies $b^2=c(s-a)$.
}
{}
{$\angle ABC=\angle ACD$.}
{\citeproblem{Suppa7049}}
{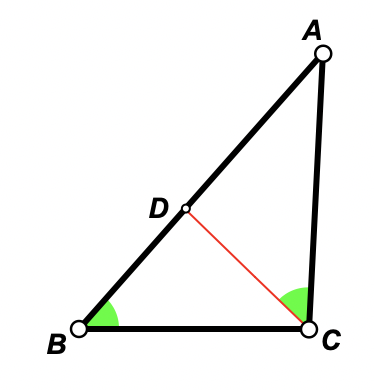}{0.4}

\bigskip
\subsection{1-step Properties}\ \medskip
\resetPropertyNumber

\prop{}
{\\Triangle $ABC$ with Gergonne cevian $BF$.\\
$\triangle ABC$ satisfies $\angle B=45\degrees$ and $\angle C=30\degrees$.\\
$D=\Gergonne[ABC]$
}
{$BE$ bisects $\angle ABC$.}
{$AF=FE$.}
{\citeproblem{RG5938}\noproof}
{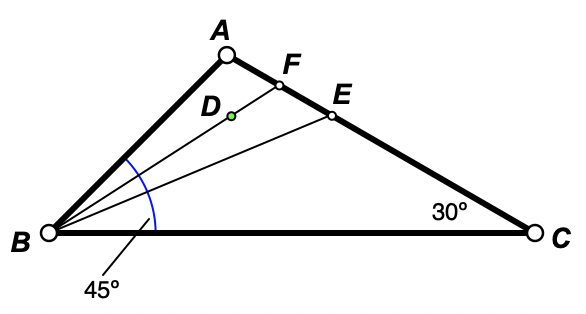}{0.5}

\prop{}
{\\Triangle $ABC$ with Gergonne cevian $AD$.
}
{The $B$-mixtilinear incircle touches $\odot ABC$ at $E$.}
{$AB\cdot CE=BE\cdot CD$.}
{\citeproblem{RG6079}}
{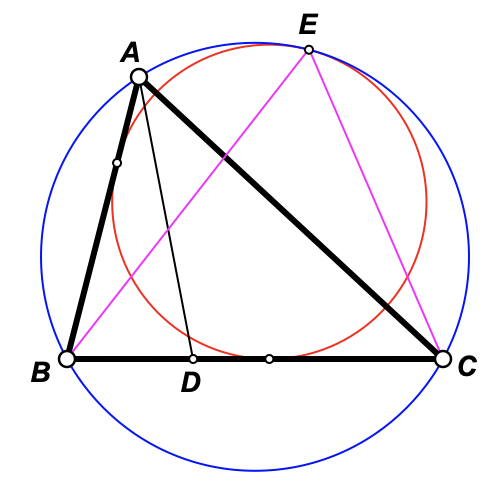}{0.4}

\newpage
\prop{}
{\\Triangle $ABC$ with Gergonne cevian $AD$.
}
{The $A$-mixtilinear incircle touches $\odot ABC$ at $E$.}
{$r(BED)=r(CDE)$.}
{\cite[Theorem~3.4]{Rabinowitz-SJM}}
{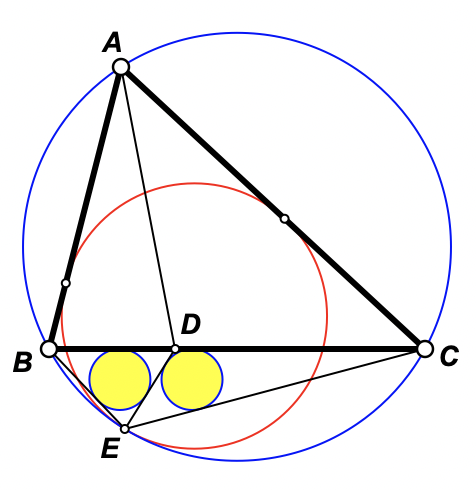}{0.4}

\prop{}
{\\Triangle $ABC$ with Gergonne cevian $AD$.
}
{The incircle of $\triangle ABD$ touches $AD$ at $E$.}
{$AC+DE=AE+CD$.}
{This follows from the split perimeter property.}
{GCIncircle}{0.4}

\prop{}
{\\Triangle $ABC$ with Gergonne cevian $CD$.\\
$\triangle ABC$ satisfies $b+c=3a$.
}
{$E$ is $A$-excenter of $\triangle ABC$.}
{$[AED]=[BCE]$.}
{\cite{vanLamoen6969}}
{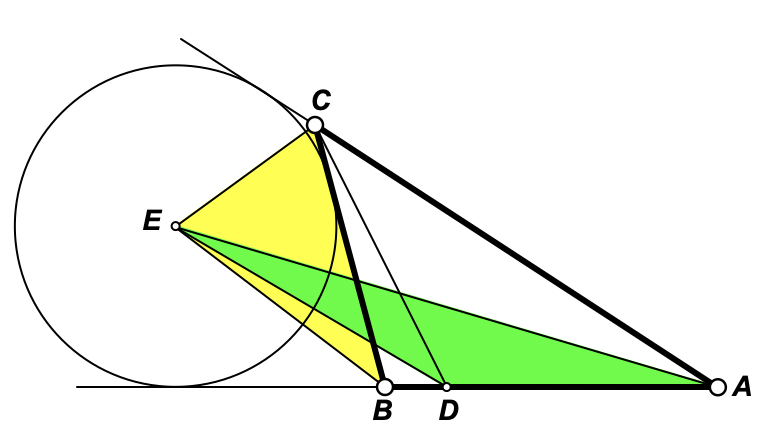}{0.4}

\prop{}
{\\Triangle $ABC$ with Gergonne cevian $AD$.\\
$\triangle ABC$ satisfies $\angle C=2\angle B$.
}
{$E=\perpen[D,AB].$}
{$\displaystyle\frac{1}{BD}+\frac{1}{BE}=\frac{1}{CD}$.}
{\citeproblem{RG6971}}
{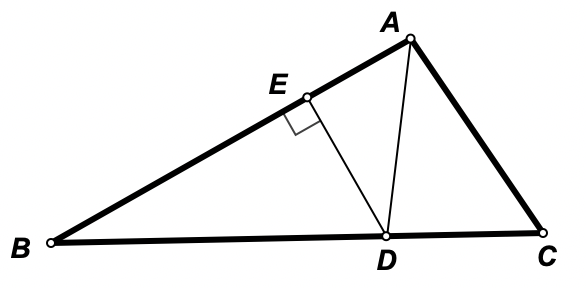}{0.5}

\newpage
\prop{}
{\\Triangle $ABC$ with Gergonne cevian $AD$.\\
$\triangle ABC$ satisfies $\angle C=2\angle B$.
}
{$E$ is the center of the $D$-excircle of $\triangle ADC$.}
{$\angle ABC=\angle AED$.}
{\citeproblem{RG6972}}
{GCDoubleAngle}{0.5}

\prop{}
{\\Triangle $ABC$ with Gergonne cevian $AD$.\\
$\triangle ABC$ satisfies $\angle B:\angle C:\angle A=1:2:4$.
}
{$E=\perpen[AC,A]\cap BC$.}
{$BD+BE=3CD$.}
{\citeproblem{RG6976}}
{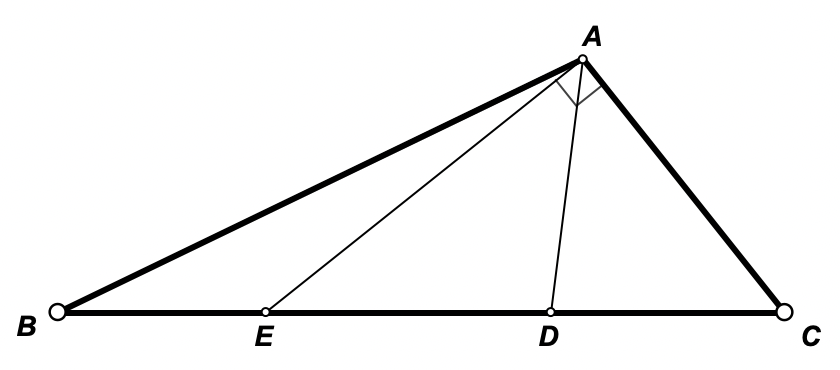}{0.5}

%******************
%
%     Triangle with Gergonne cevian and two centers
%
%******************

\bigskip
\subsection{A Triangle with a Gergonne Cevian and Two Centers}\ \medskip
\resetPropertyNumber

In this section, we report on properties found for a starting figure consisting of a
triangle and one of its Gergonne cevians along with a triangle center constructed inside
each of the two triangles formed.

\prop{}
{\\Triangle $ABC$ with Gergonne cevian $AD$.\\
$E=\incenter[ABD]$.\\
$F=\incenter[ADC]$.
}
{}
{$AD\perp EF$.}
{\citeproblem{RG6062}, follows from the Kissing Circles Theorem.
}
{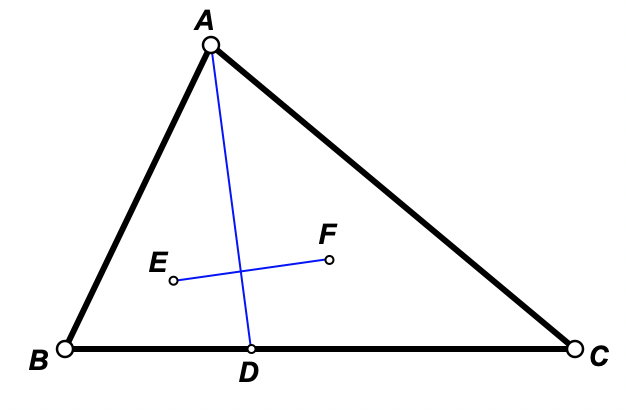}{0.45}

\prop{}
{\\Triangle $ABC$ with Gergonne cevian $AD$.\\
$E=\Gergonne[ABD]$.\\
$F=\Gergonne[ADC]$.
}
{}
{$\concur[AD,BE,CF]$.}
{This follows from the previous property.
}
{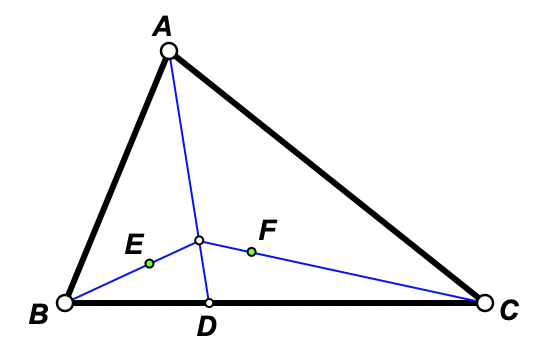}{0.4}

\prop{}
{\\Triangle $ABC$ with Gergonne cevian $AD$.\\
$E=\Nagel[ABD]$.\\
$F=\Nagel[ADC]$.
}
{}
{$\concur[AD,BE,CF]$.}
{\citeproblem{RG6063}}
{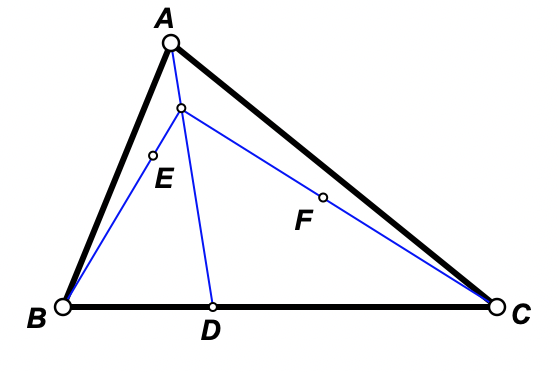}{0.4}

\prop{}
{\\Triangle $ABC$ with Gergonne cevian $AD$.\\
$\triangle ABC$ satisfies $b+c=ka$.
$E=\incenter[ABD]$.\\
$F=\centroid[ADC]$.
}
{}
{$\displaystyle\frac{[BEF]}{[BCE]}=\frac{k-3}{6}$.}
{\citeproblem{RG6995}\noproof}
{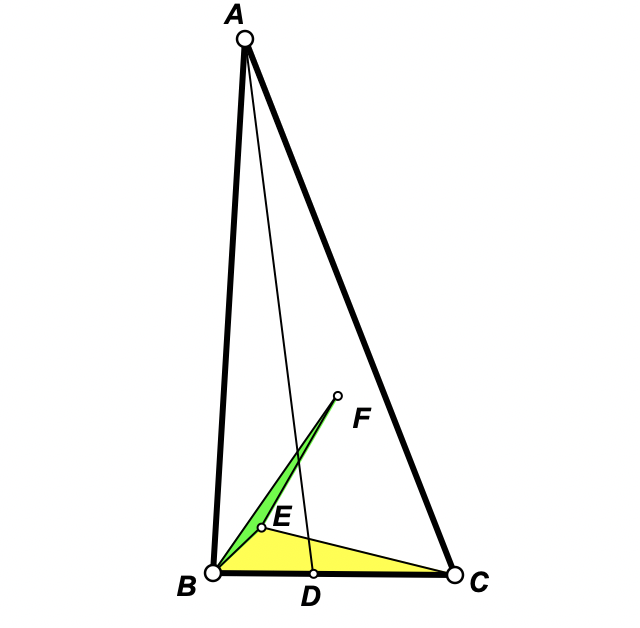}{0.4}

\prop{}
{\\Triangle $ABC$ with Gergonne cevian $AD$.\\
$\triangle ABC$ satisfies $b+c=3a$.
$E=\incenter[ABD]$.\\
$F=\centroid[ADC]$.
}
{}
{$\colline[B,E,F]$.}
{Special case of $k=3$ in previous property.}
{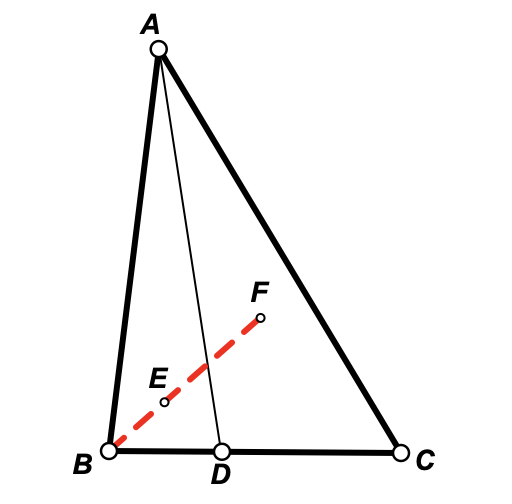}{0.4}

\newpage
\prop{}
{\\Triangle $ABC$ with Gergonne cevian $AD$.\\
$\triangle ABC$ satisfies $b+c=3a$.
$E=\centroid[ABD]$.\\
$F=\Nagel[ADC]$.
}
{}
{Three colored regions have the same area.}
{\citeproblem{RG6996}\noproof}
{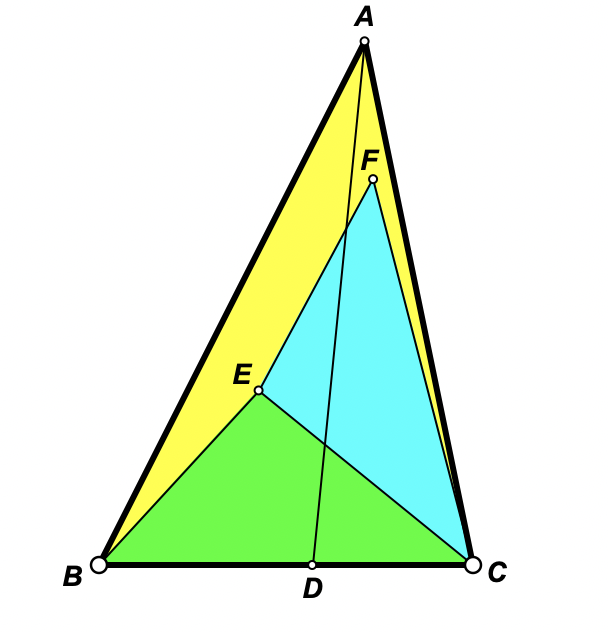}{0.4}

\prop{}
{\\Triangle $ABC$ with Gergonne cevian $AD$.\\
$\triangle ABC$ satisfies $b+c=3a$.
$E=\centroid[ABD]$.\\
$F=\Spieker[ADC]$.
}
{}
{$[BCE]=2[CFE]$.}
{\citeproblem{RG7001}}
{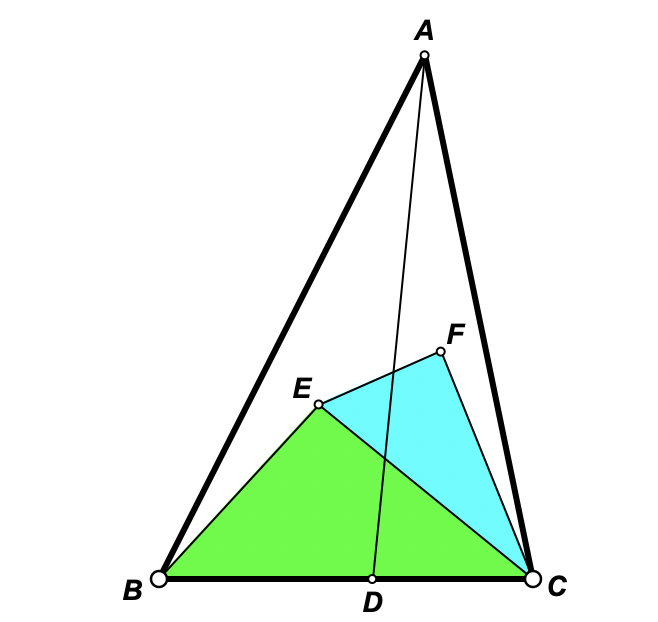}{0.4}

\prop{}
{\\Triangle $ABC$ with Gergonne cevian $AD$.\\
$\triangle ABC$ satisfies $b+c=3a$.
$E=\Gergonne[ABD]$.\\
$F=\Spieker[ADC]$.
}
{}
{$3[AEC]=8[AEF]$.}
{\citeproblem{RG7000}\noproof}
{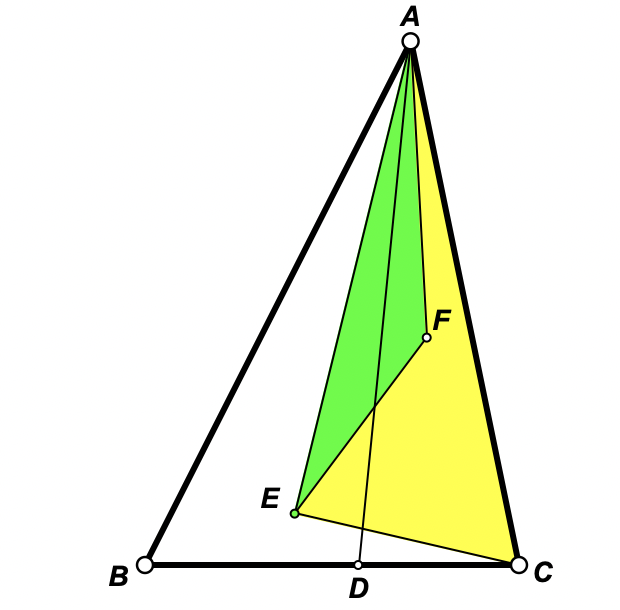}{0.4}

%******************
%
%     Triangle with Gergonne and Nagel Cevians
%
%******************

\newpage
\section{A Triangle with Gergonne and Nagel Cevians}\ \bigskip

In this section, we report on properties found for a starting figure consisting of a
triangle with two Gergonne cevians or a triangle with a Gergonne and a Nagel cevian.
 
\subsection{Two Gergonne Cevians}\ \medskip
\resetPropertyNumber

\prop{}
{Triangle $ABC$.\\
$AD$ is a Gergonne cevian.\\
$BE$ is a Gergonne cevian.
}
{}
{$CD=CE$.}
{This follows from the Trace Properties.}
{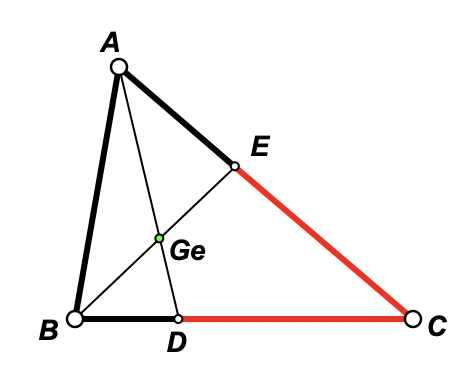}{0.4}

\prop{}
{Triangle $ABC$.\\
$\triangle ABC$ satisfies $c=2(b-a)$.\\
$BE$ is a Gergonne cevian.\\
$CD$ is a Gergonne cevian.
}
{}
{$\angle CDE=\angle BAC$.}
{\cite{Capitan7018}}
{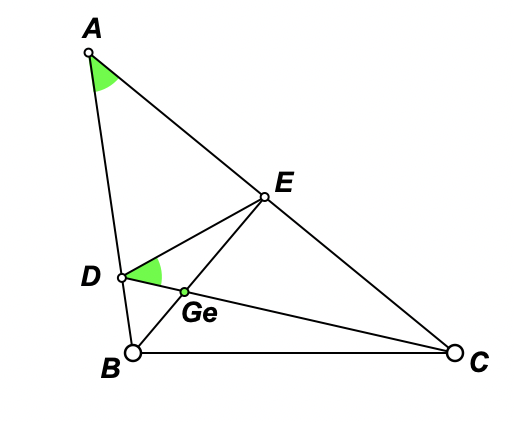}{0.4}
% also $\angle AED=\angle BDC$

\prop{}
{Triangle $ABC$.\\
$\triangle ABC$ satisfies $b(s-b)=2(s-c)(s-a)$.\\
$AE$ is a Gergonne cevian.\\
$CF$ is a Gergonne cevian.
}
{}
{$\angle ABC=\angle ADF$.}
{\citeproblem{Knop7050}}
{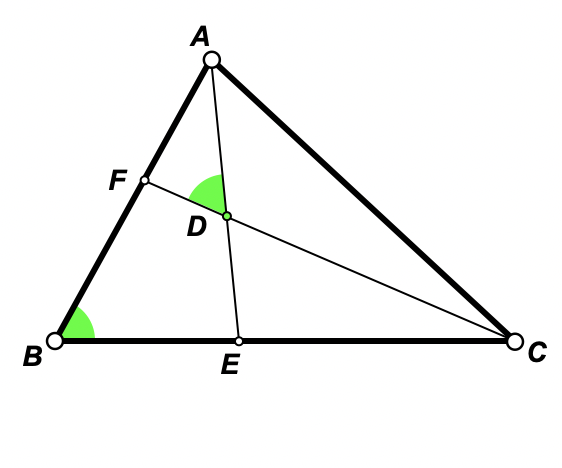}{0.4}

\prop{}
{Right triangle $ABC$.\\
$\triangle ABC$ satisfies $a=3(b-c)$.\\
$BE$ is a Gergonne cevian.\\
$CD$ is a Gergonne cevian.
}
{}
{$\odot BCD\cong\odot ADE$.}
{\cite{Suppa7018}}
{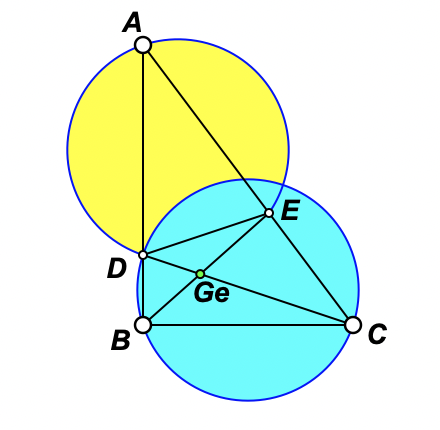}{0.4}
% Capitan has a complex generalization

\prop{}
{Triangle $ABC$.\\
$AD$ is a Gergonne cevian.\\
$BE$ is a Gergonne cevian.
}
{}
{$\displaystyle \frac{[ADC]}{[BCE]}=\frac{b}{a}$.}
{\cite{Ortega7019}}
{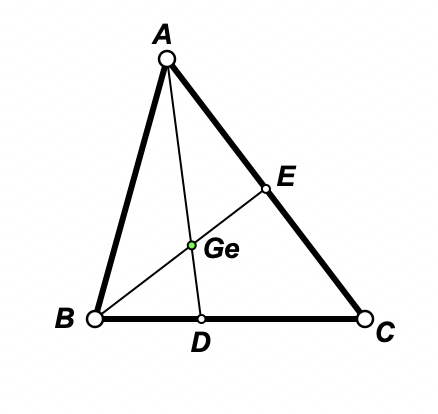}{0.4}

\prop{}
{Triangle $ABC$.\\
$\triangle ABC$ satisfies $a:b:c=16:15:11$.\\
$AD$ is a Gergonne cevian.\\
$BE$ is a Gergonne cevian.
}
{}
{$AD=CE$.\\
\vspace*{4pt}
\textbf{Note:} This result is also true whenever
$4s^3+3abc+a^2b=4s(ab+bc+ca)$.}
{\citeproblem{RG7051}}
{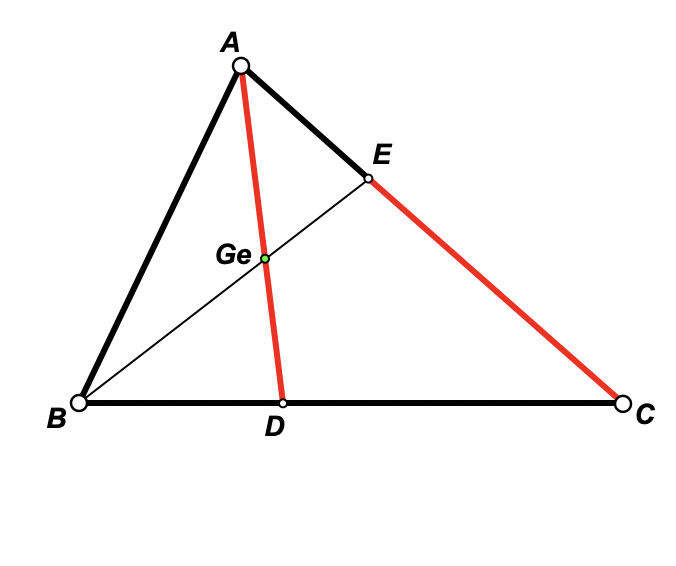}{0.4}

\newpage
\subsection{A Gergonne Cevian and a Nagel Cevian}\ \medskip
\resetPropertyNumber

\prop{Isotomic Conjugate}
{Triangle $ABC$.\\
$AD$ is a Gergonne cevian.\\
$AE$ is a Nagel cevian.
}
{}
{$BD=CE$.}
{\cite[p.~184]{Johnson}}
{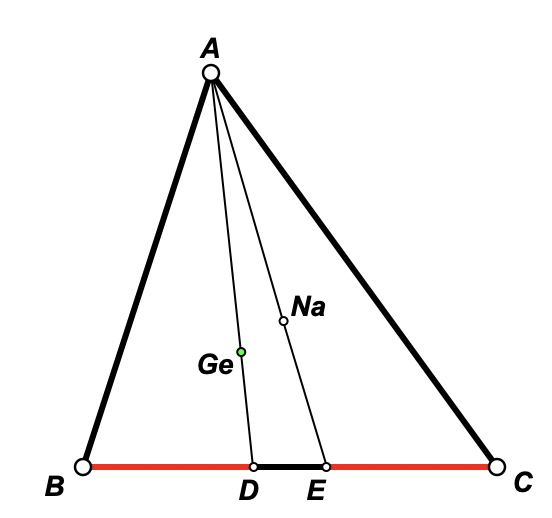}{0.4}

\prop{}
{Triangle $ABC$.\\
$AD$ is a Gergonne cevian.\\
$AE$ is a Nagel cevian.
}
{}
{$|AB-AC|=DE$.\\
Triangles $ABC$ and $ADE$ have the same centroid.}
{This follows from the Trace Properties.}
{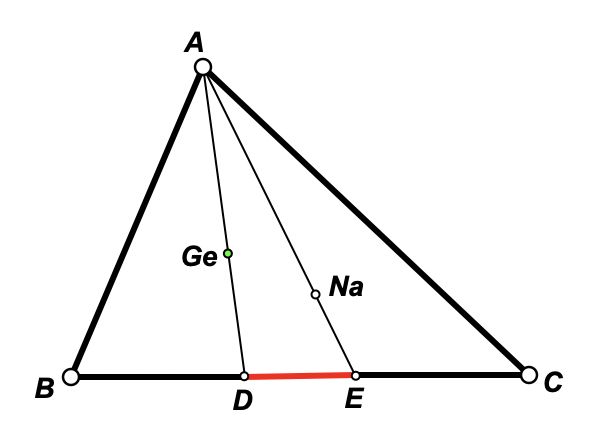}{0.4}

\prop{}
{Triangle $ABC$.\\
$AD$ is a Gergonne cevian.\\
$BE$ is a Nagel cevian.
}
{}
{$AE=CD$.}
{Follows from Trace Properties}
{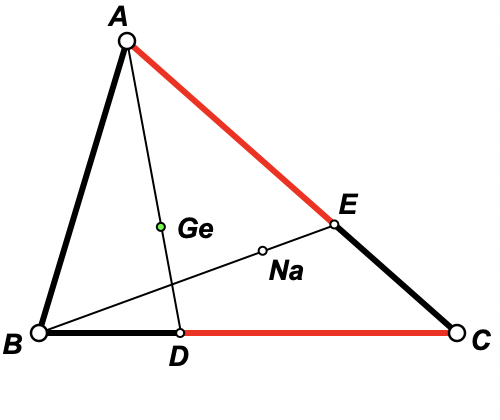}{0.4}

\newpage
\prop{}
{Right triangle $ABC$.\\
$\triangle ABC$ satisfies $\angle A=30\degrees$ and $\angle C=60\degrees$.\\
$AD$ is a Gergonne cevian.\\
$CE$ is a Nagel cevian.
}
{}
{$r(ADC)=r(BCE)$.}
{\citeproblem{RG5530}}
{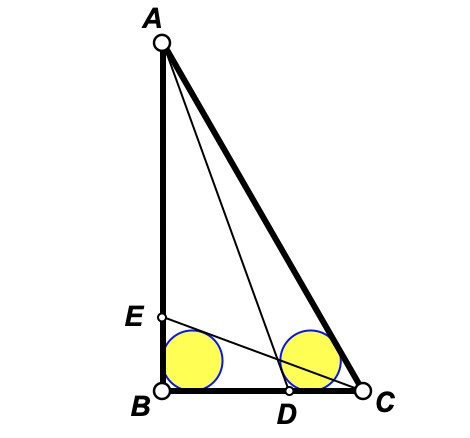}{0.4}

\prop{}
{Triangle $ABC$.\\
$AD$ is a Gergonne cevian.\\
$BE$ is a Nagel cevian.
}
{}
{$\displaystyle\frac{[ABE]}{[CED]}=\frac{2a}{b+c-a}$.}
{\cite{Capitan7017}}
{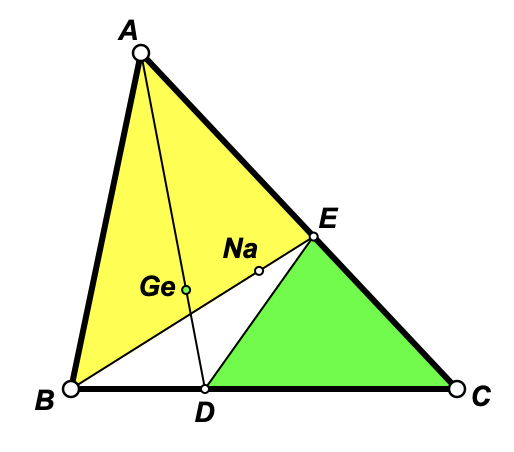}{0.4}

\prop{}
{Triangle $ABC$.\\
$AD$ is a Gergonne cevian.\\
$BE$ is a Nagel cevian.
}
{}
{$\displaystyle\frac{[ADC]}{[ABE]}=\frac{b}{a}$.}
{\cite{Marinescu7020}}
{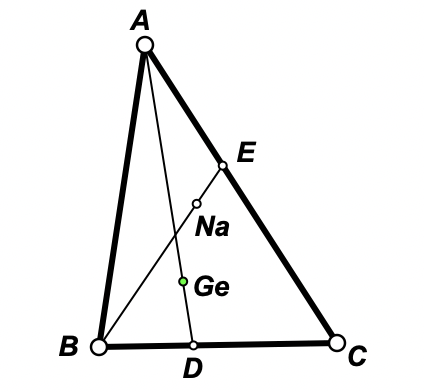}{0.4}

\prop{}
{Triangle $ABC$.\\
$\triangle ABC$ satisfies $a:b:c=5:5:2$.\\
$AD$ is a Gergonne cevian.\\
$BE$ is a Nagel cevian.
}
{}
{$BE=2AD$ and $BN_a=3AG_e$.}
{\citeproblem{RG7034}}
{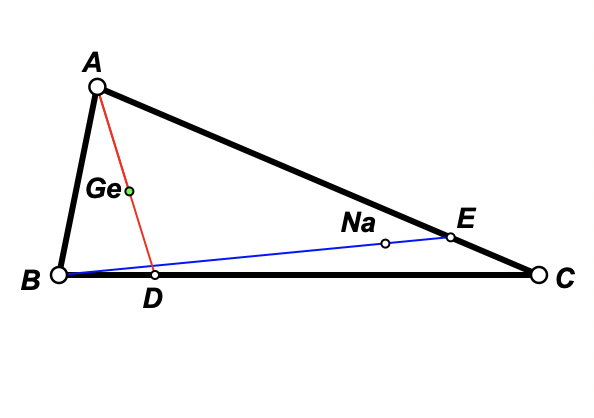}{0.4}

\prop{}
{Triangle $ABC$.\\
$\triangle ABC$ satisfies $a:b:c=3:2:2$.\\
$AD$ is a Gergonne cevian.\\
$BE$ is a Nagel cevian.
}
{}
{$BE=2AD$.\\
\vspace*{6pt}
\textbf{Note:} As in any triangle with $AB=AC$,\\
$A$, $G_e$, and $N_a$ are collinear.}
{\citeproblem{RG7035}}
{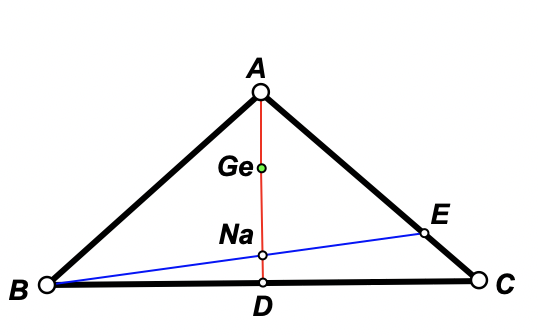}{0.5}

\prop{}
{Triangle $ABC$.\\
$\triangle ABC$ satisfies $a:b:c=9:8:3$.\\
$AD$ is a Gergonne cevian.\\
$BE$ is a Nagel cevian.
}
{}
{$BE=3AD$.\\
\vspace*{4pt}
\textbf{Note:} $BE=3AD$ is also true for triangles where $a:b:c=\sqrt{41}-3:2:2$.}
{\citeproblem{RG7036}}
{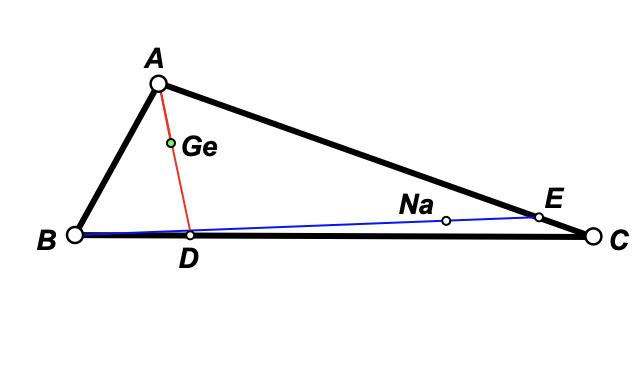}{0.45}

\prop{}
{Triangle $ABC$.\\
$\triangle ABC$ satisfies $a:b:c=9:8:5$.\\
$D=\Gergonne[ABC]$.\\
$E=\Nagel[ABC]$.
}
{}
{$AP\perp BE$ and $BQ\perp AE$.}
{\citeproblem{RG7037}}
{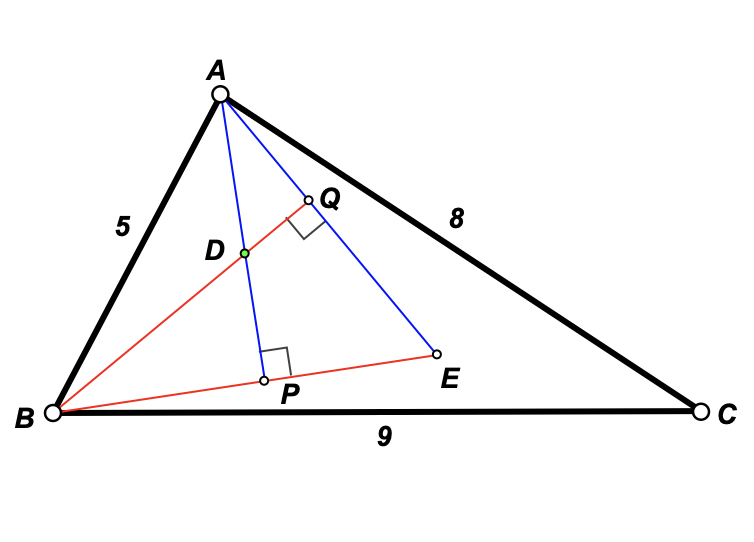}{0.4}

\prop{}
{Triangle $ABC$.\\
$D=\Gergonne[ABC]$.\\
$E=\Nagel[ABC]$.
}
{}
{$AD\perp BE$ if and only if $BD\perp AE$.}
{\citeproblem{Suppa7066}}
{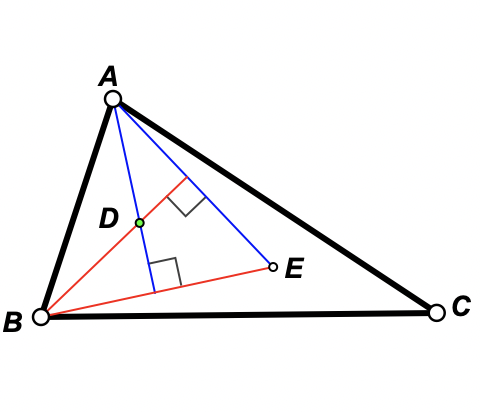}{0.4}

%******************
%
%     Triangle with cevian
%
%******************

\newpage
\section{A Triangle with a Cevian}

In this section, we report on properties found for a starting figure that consists of a
triangle $ABC$ and one of its cevians.
The cevian divides $\triangle ABC$ into two smaller triangles.
A Gergonne point is constructed in one of these triangles or in the original triangle.

\subsection{Gergonne Point in Original Triangle}\ \medskip
\resetPropertyNumber

\prop{}
{Right triangle $ABC$ with cevian through the nine-point center.\\
$\triangle ABC$ satisfies $c=2(b-a)$ and $\angle B=90\degrees$.\\
$D=\Gergonne[ABC]$.\\
$E=\ninepointcenter[ABC]$.
}
{}
{\\
(a) $\colline[A,E,D]$.\\
(b) $\displaystyle\frac{AE}{ED}=\frac{3(2a+b)}{9b-14a}$.}
{\cite{Suppa6007}}
{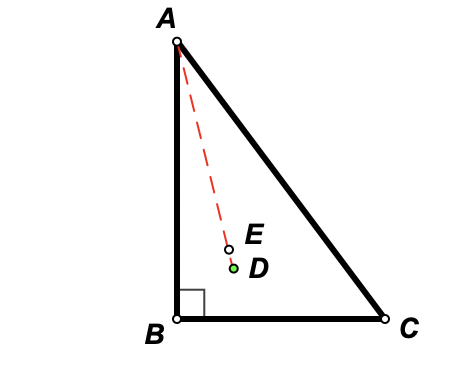}{0.43}

\prop{}
{Triangle $ABC$ with median.\\
$\triangle ABC$ satisfies $c=2(b-a)$\\
$D=\Gergonne[ABC]$.\\
$E=\midpoint[AB]$.
}
{}
{$[BCD]=[CED]$.}
{\citeproblem{Capitan6009}\noproof}
{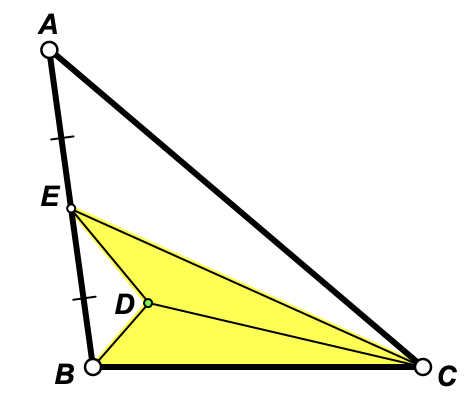}{0.4}

\prop{}
{Triangle $ABC$ with median.\\
$\triangle ABC$ satisfies $b+c=2a$.\\
$D=\Gergonne[ABC]$.\\
$E=\midpoint[BC]$.
}
{}
{$(AD)^2+(AE)^2=(DE)^2+(BC)^2$.}
{\citeproblem{RG6006}\noproof}
{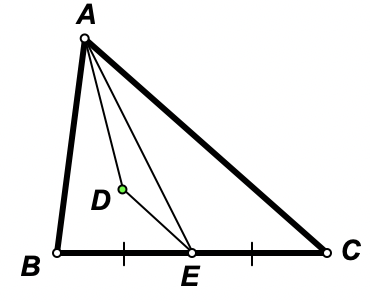}{0.4}

\prop{Isogonal Conjugate}
{\\Triangle $ABC$ with an $X_{55}$-cevian.\\
$G_e=\Gergonne[ABC]$.
}
{}
{$\angle BAG_e=\angle CAX_{55}$\\
\vspace*{8pt}
\textbf{Note:} $I$ is the incenter, $O$ is the circumcenter, $M$ is the centroid, and $F_e$ is the Feuerbach point of $\triangle ABC$.}
{\cite{ETC55}}
{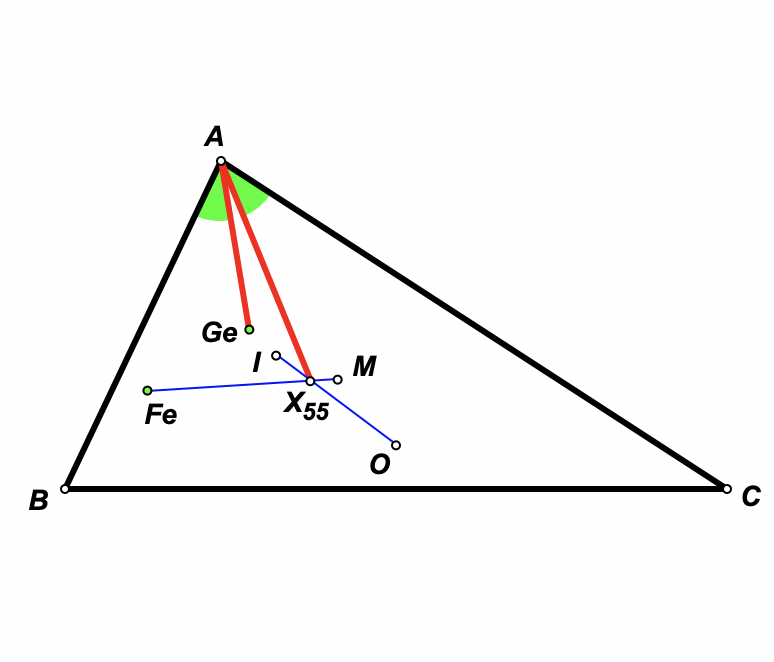}{0.4}

\bigskip
\subsection{Gergonne Point in One of the Smaller Triangles}\ \medskip
\resetPropertyNumber

\prop{}
{\\Triangle $ABC$ with Nagel cevian $BD$.\\
$\triangle ABC$ satisfies $b+c=3a$.\\
$E=\Gergonne[BCD]$.
}
{}
{$CE\perp BD$.}
{\citeproblem{RG7002}}
{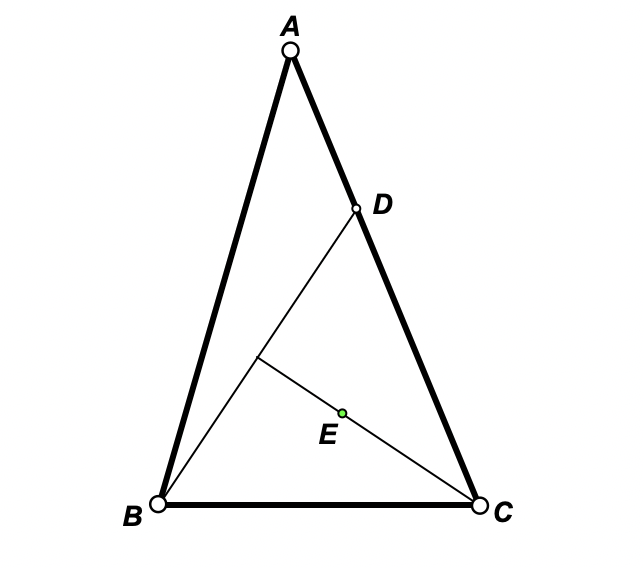}{0.4}

\prop{}
{\\Triangle $ABC$ with symmedian $CD$.\\
$\triangle ABC$ satisfies $a^2+b^2=bc$.\\
$E=\Gergonne[ADC]$.
}
{}
{$AE\perp CD$.}
{\cite{Suppa7003}}
{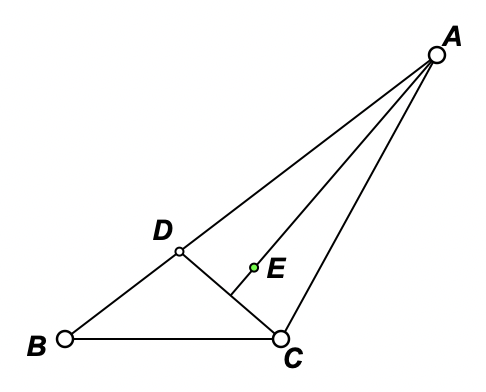}{0.4}

\newpage
\prop{}
{\\Triangle $ABC$ with $X_3$-cevian $AD$.\\
$\triangle ABC$ satisfies $\angle C=2\angle B$.\\
$E=\Gergonne[ADC]$.
}
{}
{$CE\perp AD$.}
{\citeproblem{RG7004}}
{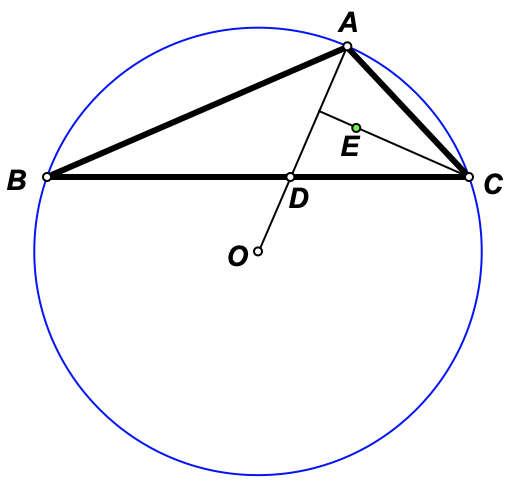}{0.4}

\prop{}
{\\Right triangle $ABC$ with angle bisector $BD$.\\
$\triangle ABC$ satisfies $a:b:c=3:4:5$.\\
$E=\Gergonne[BCD]$.
}
{}
{\\
(a) $\angle BEC=\angle ADB$.\\
(b) $(DE)^2+(BD)^2=(AE)^2$.\noproof
}
{\citeproblem{RG6069}}
{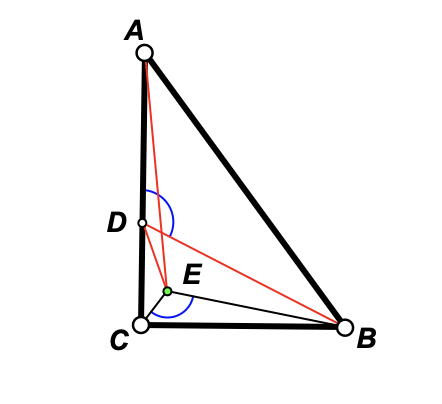}{0.4}

\prop{}
{\\Triangle $ABC$ with angle bisector $AD$.\\
$\triangle ABC$ satisfies $a:b:c=6:5:4$.\\
$E=\Gergonne[ACD]$.
}
{}
{\\
(a) $\angle AEC$ and $\angle ABC$ are supplementary.\noproof\\
(b) $(BE)^2+(AC)^2=(BC)^2$.\noproof\\
(c) $E=\Nagel[ABC]$.
}
{\citeproblem{RG6070}}
{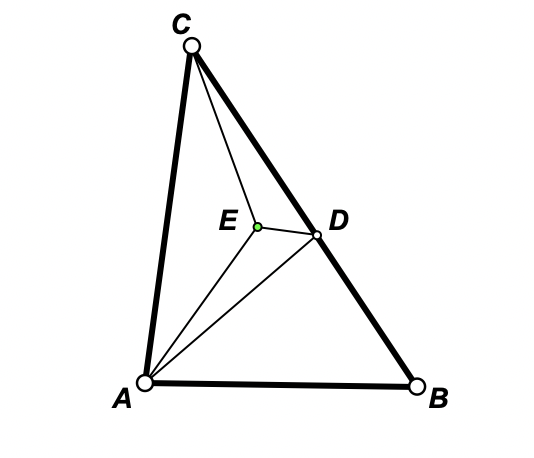}{0.4}

\prop{}
{\\Triangle $ABC$ with angle bisector $AD$.\\
$\triangle ABC$ satisfies $a:b:c=7:5:3$.\\
$E=\Gergonne[ACD]$.
}
{}
{$[ABD]=[AEC]$.}
{\citeproblem{RG6072}}
{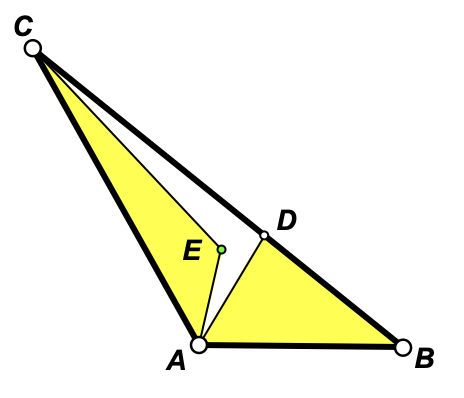}{0.4}

\prop{}
{\\Triangle $ABC$ with angle bisector $AD$.\\
$\triangle ABC$ satisfies $a:b:c=7:6:3$.\\
$E=\Gergonne[ACD]$.
}
{}
{\\
(a) $AB\parallel DE$.\hfill
(b) $BD+DE=AB$.\\
(c) $AC=9DE$.\hfill
(d) $CD=7DE$.\\
(e) $CD+DE=2AD$.\hfill
(f) $AD=4DE$.
}
{\citeproblem{RG6071}\noproof}
{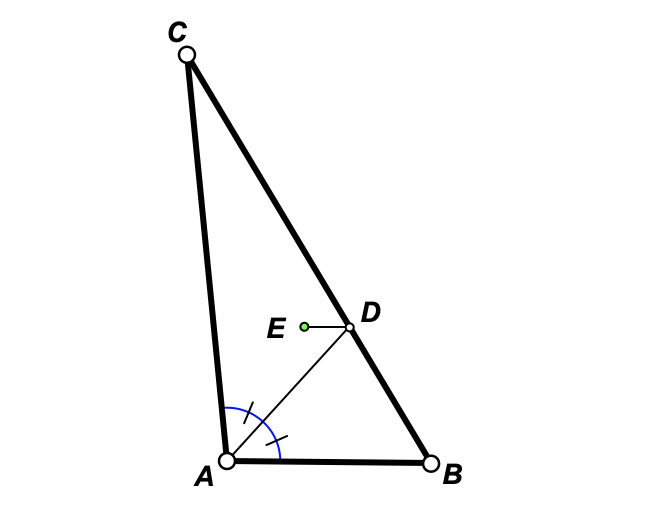}{0.4}

\prop{}
{\\Triangle $ABC$ with angle bisector $CD$.\\
$\triangle ABC$ satisfies $a:b:c=9:6:5$.\\
$E=\Gergonne[ACD]$.
}
{}
{$\angle ABE=\angle ACE$.}
{\citeproblem{RG6073}}
{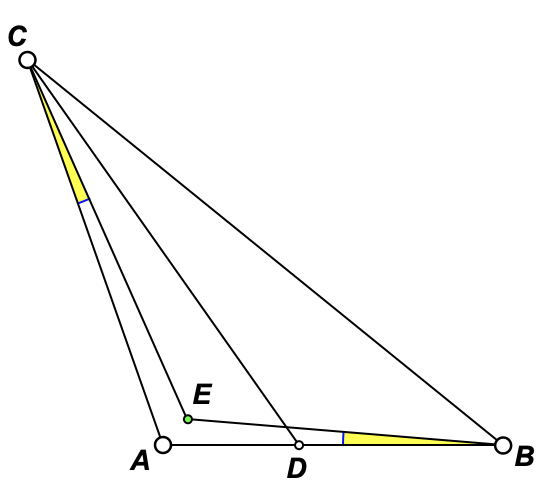}{0.4}

\prop{}
{\\Triangle $ABC$ with angle bisector $AD$.\\
$\triangle ABC$ satisfies $a^2=c(b+c)$.\\
$E=\Gergonne[ADC]$.
}
{}
{$DE$ is the perpendicular bisector of $AC$.
}
{\citeproblem{RG6070}}
{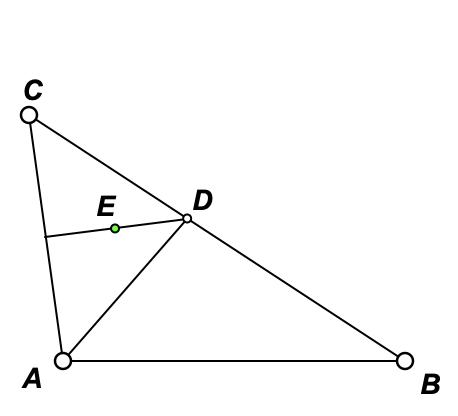}{0.4}

%******************
%
%     Triangle with Gergonne Pararadius
%
%******************

\newpage
\section{Lines Through the Gergonne Point}

\begin{figure}[h!t]
\centering
\includegraphics[width=0.4\linewidth]{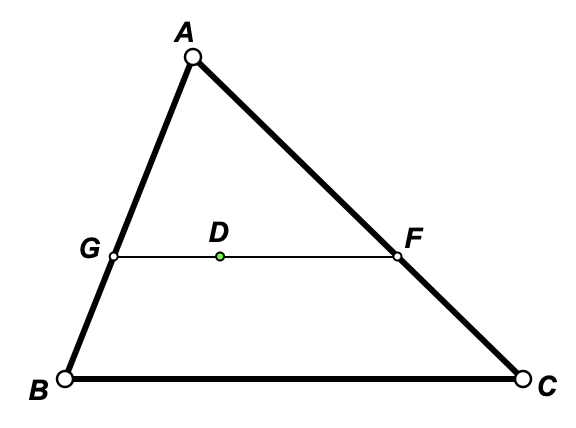}
\includegraphics[width=0.4\linewidth]{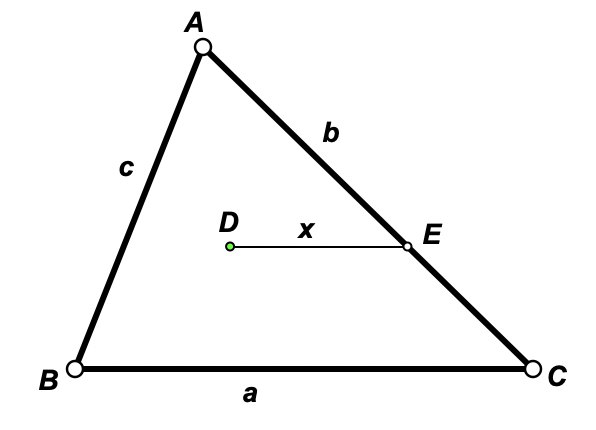}
\caption{Gergonne Chord and Gergonne Pararadius}
\label{fig:parachord}
\end{figure}

The line segment joining two points on the perimeter of a triangle
is called a \textit{chord}.
A chord parallel to a side of the triangle is called a \textit{parachord}.
This is also known as a \textit{parallelian} in the literature.
The parachord through the Gergonne point of a triangle is
called a \textit{Gergonne parachord} (Figure~\ref{fig:parachord}, left). The line segment from the Gergonne point
to the point where the Gergonne parachord meets a side of the triangle is called 
a Gergonne pararadius (Figure~\ref{fig:parachord}, right).

\subsection{Properties of a Pararadius}\ \medskip
\resetPropertyNumber

In this section, we report on properties found for a starting figure consisting of a
triangle and a Gergonne pararadius.

\prop{(Pararadius Formula)}
{\\Triangle $ABC$ with Gergonne pararadius.\\
$D=\Gergonne[ABC]$.\\
$E=\para[D,BC]\cap AC$.
}
{}
{\\
\vspace*{4pt}
$\displaystyle x=\frac{-a(a+b-c)(b+c-a)}{a^2+b^2+c^2-2ab-2bc-2ca}$.}
{See supplementary notebook.}
{pararadius}{0.47}

\prop{(Parallel Ratio)}
{\\Triangle $ABC$ with Gergonne pararadius.\\
$D=\Gergonne[ABC]$.\\
$E=\para[D,BC]\cap AC$.
}
{}
{$\displaystyle \frac{AE}{DE}=\frac{b}{s-c}$.}
{See supplementary notebook.}
{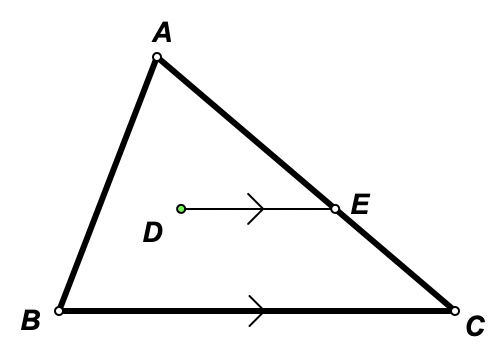}{0.4}

\prop{}
{\\Triangle $ABC$ with Gergonne pararadius.\\
$D=\Gergonne[ABC]$.\\
$E=\para[D,BC]\cap AC$.
}
{}
{\\
\vspace*{4pt}
$\displaystyle AE=\frac{2ab(a-b-c)}{a^2+b^2+c^2-2ab-2bc-2ca}$.\\
\vspace*{4pt}
$\displaystyle CE=\frac{b(b^2+c^2-a^2-2bc)}{a^2+b^2+c^2-2ab-2bc-2ca}$.}
{See supplementary notebook.}
{pararadius}{0.6}

\prop{}
{\\Triangle $ABC$ with Gergonne pararadius.\\
$\triangle ABC$ satisfies $\angle B=90\degrees$.\\
$D=\Gergonne[ABC]$.\\
$E=\para[D,BC]\cap AC$.}
{}
{$\displaystyle BC+AE=\left(1+\frac{2b}{a}\right)DE$.}
{\cite{Suppa7052}. See also \citeproblem{RG5937}.}
{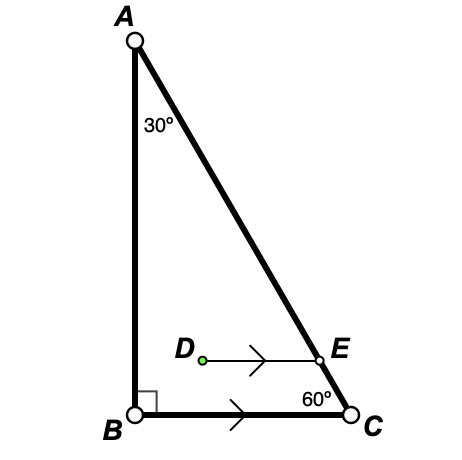}{0.4}

\prop{}
{\\Triangle $ABC$ with Gergonne pararadius.\\
$\triangle ABC$ satisfies $c=3(b-a)$.\\
$D=\Gergonne[ABC]$.\\
$E=\para[D,BC]\cap AC$.}
{}
{$\displaystyle \frac{1}{AE}+\frac{1}{DE}=\frac{1}{CE}$.}
{\citeproblem{Suppa6010}}
{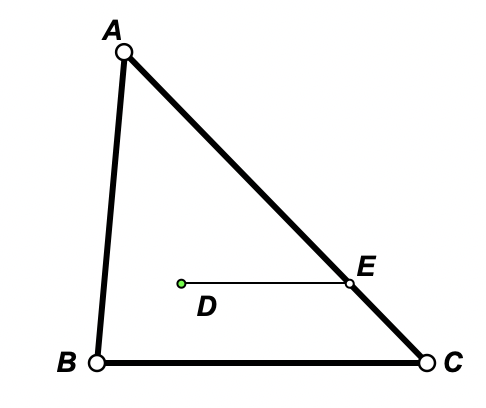}{0.4}

\prop{}
{\\Triangle $ABC$ with Gergonne pararadius.\\
$\triangle ABC$ satisfies $a:b:c=3:5:4$.\\
$D=\Gergonne[ABC]$.\\
$E=\para[D,AB]\cap AC$.}
{}
{The circle through $D$ touching $AC$ at $E$ is tangent to $BC$.}
{\citeproblem{RG6011}}
{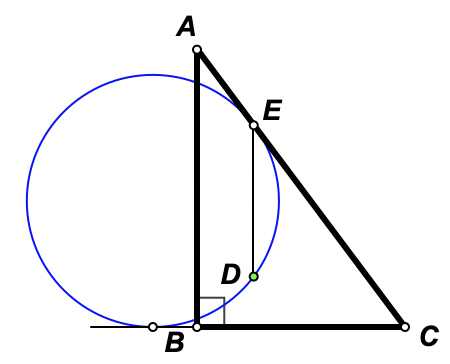}{0.4}

\prop{}
{\\Triangle $ABC$ with Gergonne pararadius.\\
$\triangle ABC$ satisfies $c=(2-\sqrt5)a+b$.\\
$D=\Gergonne[ABC]$.\\
$E=\para[D,AC]\cap BC$.
}
{}
{
$\displaystyle x=\frac{2b(c-b)}{a-3b+c}$.}
{\citeproblem{RG6017} and \cite{Suppa5982}}
{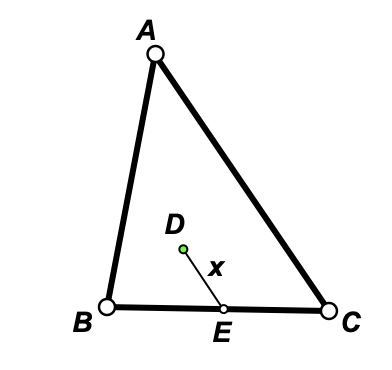}{0.47}

\prop{}
{\\Triangle $ABC$ with three Gergonne pararadii.\\
$P=\Gergonne[ABC]$.\\
$X=\para[P,AB]\cap BC$.\\
$Y=\para[P,BC]\cap CA$.\\
$Z=\para[P,CA]\cap AB$.
}
{}
{
\vspace*{8pt}
\\$\displaystyle \frac{1}{AZ}+\frac{1}{PX}=\frac{1}{BX}+\frac{1}{PY}=\frac{1}{CY}+\frac{1}{PZ}$.}
{\citeproblem{RG6176}}
{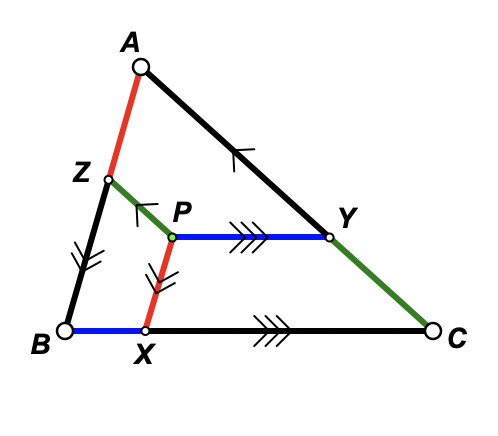}{0.43}

%******************
%
%     Triangle with Gergonne Parachord
%
%******************

\newpage
\subsection{Properties of a Parachord}\ \medskip
\resetPropertyNumber

In this section, we report on properties found for a starting figure consisting of a
triangle and a Gergonne parachord.

\prop{(Parachord Formula)}
{\\Triangle $ABC$ with Gergonne parachord.\\
$D=\Gergonne[ABC]$.\\
$F=\para[D,BC]\cap AC$.\\
$G=\para[D,BC]\cap AB$.
}
{}
{\\
\vspace*{4pt}
$\displaystyle FG=\frac{2a^2(a-b-c)}{a^2+b^2+c^2-2ab-2bc-2ca}$.}
{This follows from the Pararadius Formula.}
{parachord}{0.5}

\prop{(Parachord Perimeter)}
{\\Triangle $ABC$ with Gergonne parachord.\\
$D=\Gergonne[ABC]$.\\
$F=\para[D,BC]\cap AC$.\\
$G=\para[D,BC]\cap AB$.
}
{}
{
$AG+DF=AF+DG$.}
{This follows from the Split Perimeter Property.}
{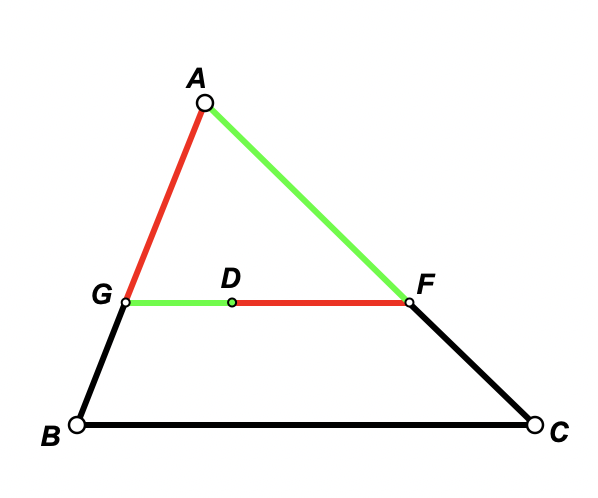}{0.4}

%******************
%
%     Triangle with Gergonne Chord
%
%******************

\bigskip
\subsection{Properties of Gergonne Chords}\ \medskip
\resetPropertyNumber

In this section, we report on properties found for a starting figure consisting of a
triangle and a Gergonne chord.

\prop{}
{Isosceles triangle $ABC$.\\
$EF$ is a Gergonne chord.\\
$\triangle ABC$ satisfies $b=c$.\\
$D=\Gergonne[ABC]$.\\
$F=ED\cap BC$.
}
{}
{\\
\vspace*{4pt}
$m(2m+2n-p-q)(p-q)=np(p+q)$.}
{\citeproblem{RG6828}\noproof}
{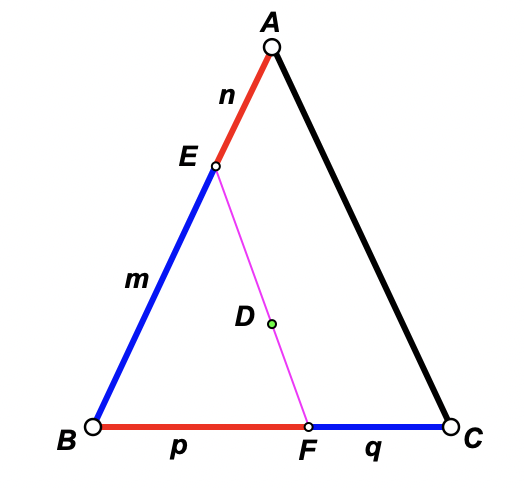}{0.4}

\prop{}
{Triangle $ABC$.\\
$EF$ is a Gergonne chord.\\
$D=\Gergonne[ABC]$.\\
$F=ED\cap AC$.
}
{}
{\\
\vspace*{8pt}
$\displaystyle\frac{s-a}{s-b}\cdot\frac{n}{m}+\frac{s-a}{s-c}\cdot\frac{q}{p}=1$.}
{\cite[p.~46]{Rosado}}
{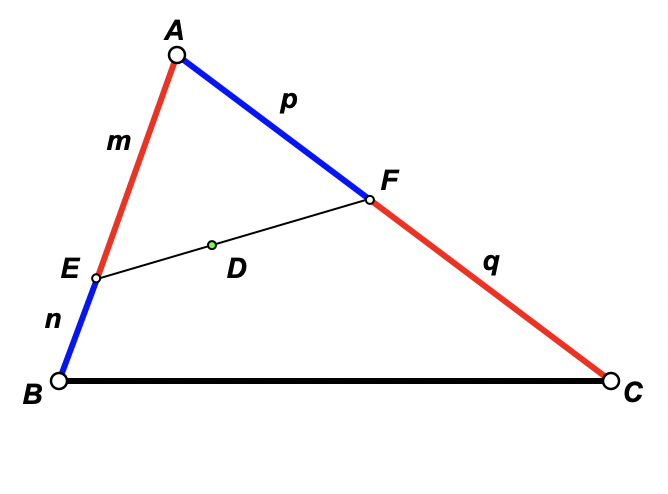}{0.5}

%******************
%
%     Perpendiculars From the Gergonne Point
%
%******************

\bigskip
\subsection{Perpendiculars From the Gergonne Point}\ \medskip
\resetPropertyNumber

In this section, we report on properties found for a starting figure consisting of a
triangle and one or more perpendiculars dropped from the Gergonne point to the sides
of the triangle.

\prop{(Apothem Length)}
{Triangle $ABC$.\\
$D=\Gergonne[ABC]$.\\
$E=\perpen[D,BC]$
}
{}
{\\
\vspace*{14pt}
$\displaystyle DE=\frac{8(s-b)(s-c)K}{a(2ab+2bc+2ca-a^2-b^2-c^2)}$,\\
\vspace*{14pt}
where $K=\sqrt{s(s-a)(s-b)(s-c)}$.
}
{This follows from the area formula.}
{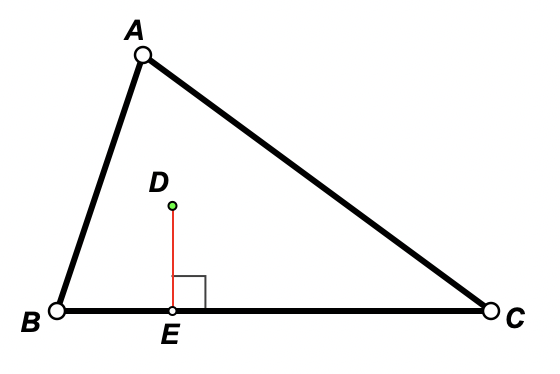}{0.6}

\prop{(Trilinear Coordinates)}
{Triangle $ABC$.\\
$D=\Gergonne[ABC]$.\\
$E=\perpen[D,BC]$\\
$F=\perpen[D,AC]$
}
{}
{
$\displaystyle \frac{DE}{DF}=\frac{b(c+a-b)}{a(b+c-a)}$.\\
}
{\cite{ETC7}}
{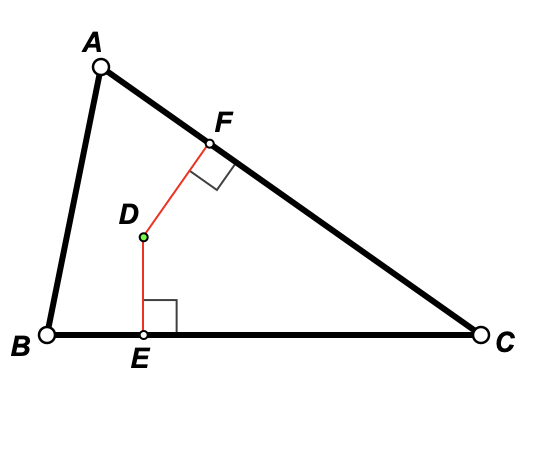}{0.4}
% This follows from the Apothem Length formula.

\prop{}
{Triangle $ABC$.\\
$\triangle ABC$ satisfies $a:b:c=16:15:7$.\\
$D=\Gergonne[ABC]$.\\
$E=\perpen[D,BC]$
}
{}
{
$CD=2AE$.\\
}
{\citeproblem{RG7054}}
{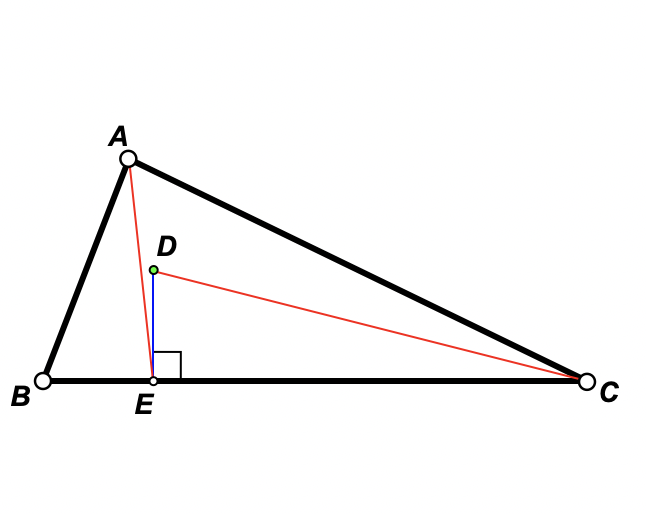}{0.4}

\prop{}
{Triangle $ABC$.\\
$\triangle ABC$ satisfies $a:b:c=8:5:4$.\\
$D=\Gergonne[ABC]$.\\
$E=\perpen[D,BC]$
}
{}
{
$\displaystyle \frac{[BCD]}{[AED]}=632$.\\
}
{\citeproblem{RG7030}}
{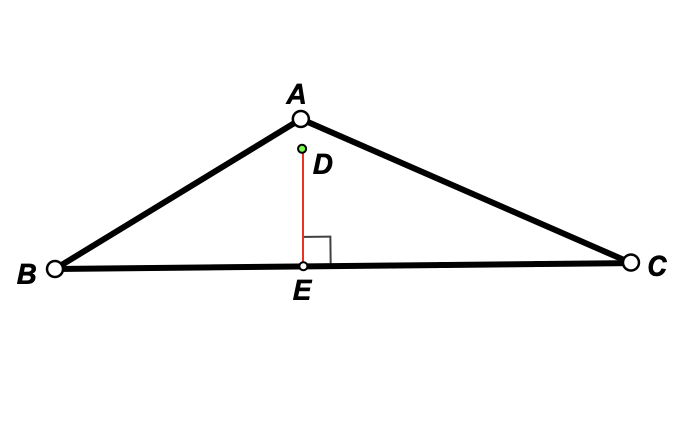}{0.4}

\prop{}
{Right triangle $ABC$.\\
$\triangle ABC$ satisfies $\angle A=90\degrees$.\\
$D=\Gergonne[ABC]$.
$E=\perpen[D,AC]$
$F=\perpen[D,AB]$
$G=\perpen[D,BC]$
}
{}
{
$\displaystyle \frac{1}{x}+\frac{1}{y}=\frac{1}{z}$.
}
{\citeproblem{RG7033}}
{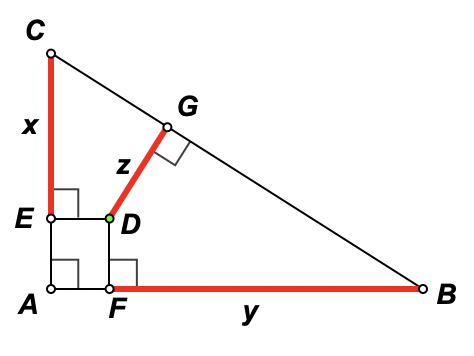}{0.4}

\bigskip
\bigskip
\bigskip
\section{Other Properties Discovered by Computer}

Other properties of the Gergonne point (discovered by computer) were found by Dekov \cite{Dekov}. A typical result is: The Gergonne point of a triangle is the mittenpunkt of the orthic triangle of the intouch triangle.

%\newpage
%==========================
% Bibliography
%==========================
%\bigskip

\goodbreak

% These involve the Gergonne point but have too many steps:
%  RG5398, RG4286, RG4691, RG4271, RG4260, RG4514,
% RG4058, RG4459, RG4463, RG4623, RG4701, RG5398, RG5454
% RG5790, RG5794, RG5800, RG4692, RG5821}
% These involve cevasix or parasix configuration
% RG3318, RG4383, RG4445, RG5269

\end{document}